\documentclass[twoside,a4paper]{amsart}
\usepackage[centertags]{amsmath}
\usepackage{amsfonts}
\usepackage{dsfont}
\usepackage{amssymb}
\usepackage{amsthm}
\usepackage[ansinew]{inputenc}
\usepackage[cmtip,all,poly]{xy}
\usepackage{graphicx}
\newtheorem{thm}[equation]{Theorem}

\newtheorem{lem}[equation]{Lemma}
\newtheorem{prop}[equation]{Proposition}
\theoremstyle{definition}
\newtheorem{defn}[equation]{Definition}
\theoremstyle{remark}
\newtheorem{rem}[equation]{Remark}

\numberwithin{equation}{section}
\newcommand{\abs}[1]{\left\vert#1\right\vert}
\newcommand{\set}[1]{\left\{#1\right\}}

\newcommand{\To}{\longrightarrow}

\newcommand{\C}[1]{\mathbf{#1}} % categor�s

\def\r{\rightarrow} % flecha -->
\def\rr{\Rightarrow} % flecha ==>

\newcommand{\sym}[1]{\operatorname{Sym}(#1)}

\newcommand{\symt}[1]{\operatorname{Sym}_\vc(#1)}

\def\sign{\operatorname{sign}}

\def\st{\stackrel} % abreviatura de \stackrel

\def\coker{\operatorname{Coker}}
\renewcommand{\ker}{\operatorname{Ker}}

\def\Z{\mathbb{Z}}

\def\sq{Sq}

\def\N{\mathds{N}}

\def\S{\Sigma}

\newcommand{\grupo}[1]{\langle #1\rangle}

\newcommand{\vc}{\Box}    % composici� vertical de 2-celdas
\begin{document}

\title{Toda brackets and cup-one squares for ring spectra}%
\author{Hans-Joachim Baues and Fernando Muro}%

\address{Max-Planck-Institut f\"ur Mathematik, Vivatsgasse 7, 53111 Bonn, Germany}%
\email{baues@mpim-bonn.mpg.de, muro@mpim-bonn.mpg.de}%

\thanks{The second author was partially supported
by the project MTM2004-01865 and the MEC postdoctoral fellowship EX2004-0616.}%

\date{\today}
%\dedicatory{}%
%\commby{}%
% ----------------------------------------------------------------	
\begin{abstract}
In this paper we prove the laws of Toda brackets on the homotopy groups of a connective ring spectrum and  the laws of the
cup-one square in the homotopy groups of a commutative connective ring spectrum.
\end{abstract}
\maketitle
\begin{footnotesize}
\tableofcontents
\end{footnotesize}
% ----------------------------------------------------------------

\section*{Introduction}

Secondary homotopy operations such as triple Toda brackets are defined on the homotopy groups of
a ring spectrum $R$ enriching the ring structure of $\pi_*R$. Toda established in \cite{toda} a set of relations for Toda brackets in
the stable homotopy groups of spheres, Alexander claimed these relations for some cobordism rings in \cite{cmp},
and we show here that the Toda relations are, in fact, satisfied for any connective ring spectrum $R$ (Definition
\ref{rmp} and Theorem \ref{main1}). Moreover, if the ring spectrum $R$ is commutative further relations proved
by Toda for the sphere spectrum, such as the Jacobi identity, are shown to be satisfied in general (Definition \ref{crmp} and Theorem \ref{main3}).

If $R$ is commutative a new secondary
homotopy operation appears, namely the cup-one square. This operation was studied in \cite{hirs} in the context
of $H_\infty$-ring spectra. The operation in \cite{hirs} is, however, only defined up to an indeterminacy. 
We show that one can extract from this undetermined operation a fully determined cup-one square and we compute its behaviour
with respect to sums and products in $\pi_*R$, as well as its relation to Toda brackets (Definition \ref{crmp} and Theorem \ref{main3}). This is done 
by carrying out a careful analysis in the ``symmetric track groups'' introduced in \cite{2hg2}. 
In this way we are able to compute explicitly the deviation of the cup-one square from additivity and from being a quadratic
derivation, which was only computed in \cite{hirs} up to an unknown constant, see (T11) and (T12) in Definition
\ref{crmp}.

For the proofs we use the algebraic framework of ($E_\infty$-)quadratic pair algebras, which are algebraic models
of (commutative) ring spectra extending the homotopy groups and codifying all secondary operations, see
\cite{2hg4}. 

The statements of the homotopical results are in the first section, which can be regarded as a continuation
of this introduction. The rest of sections are purely algebraic and contain all proofs.

\section{Secondary operations and their laws}

A \emph{(commutative) ring spectrum} is a (commutative) monoid in 
the closed symmetric monoidal model category of symmetric spectra of
compactly generated topological spaces defined in \cite[12]{mcds}. The monoidal structure is given by the 
smash product $X\wedge
Y$ and the unit object is the \emph{sphere spectrum} $S$. A symmetric spectrum is \emph{connective} if its
homotopy groups vanish in negative dimensions.

The homotopy groups of a connetive ring spectrum $\pi_*R$ form an $\N$-graded ring, where
$\N=\set{0,1,2,\dots}$. All rings and modules in this paper will be $\N$-graded and the degree of a homogeneous
element $x$ will be denoted by $\abs{x}$. Ungraded objects are regarded
as graded objects concentrated in degree $0$. The degree of a homogeneous element $a\in\pi_*R$ is denoted by $\abs{a}$. 
The ring $\pi_*R$ is equipped with secondary homotopy operations called
\emph{Toda brackets}. The Toda bracket of three homogeneous elements
$$\grupo{a,b,c}\subset \pi_{\abs{a}+\abs{b}+\abs{c}+1}R$$
is a coset of
$$(\pi_{\abs{a}+\abs{b}+1}R)\cdot c+a\cdot (\pi_{\abs{b}+\abs{c}+1}R)$$
which is defined whenever $a b=0$ and $b c = 0$.
This operation was first considered by Toda for the sphere spectrum $S$, see
\cite{toda}. In \cite{cmp} one finds a construction of Toda brackets for various cobordism spectra under the name of Massey products. 
We consider in \cite{2hg4} two equivalent definitions of Toda brackets on the homotopy groups of a ring spectrum $R$.
Both definitions use the model category of \emph{right $R$-modules}, see \cite[12]{mcds}. One of the definitions uses Toda brackets for
triangulated categories in the sense of \cite{shc} applied to the homotopy category of $R$-modules. This is also
the definition adopted in \cite{steffen}. The alternative definition uses \emph{tracks}, i.e. homotopy
classes of homotopies, in the model category of $R$-modules. We now recall this definition.

The homotopy group $\pi_nR$ coincides with the group of morphisms from the $n$-fold suspension 
$\S^nR\r R$ in the homotopy category of right $R$-modules.
We can suppose without
loss of generality that $R$ is a fibrant ring spectrum. In that case the elements $a,b,c\in\pi_*R$ can be
realized by maps $\bar{a},\bar{b},\bar{c}$ in the category of $R$-modules.
The vanishing hypothesis $ab=0$ and $bc=0$ imply the existence of null-homotopies
\begin{equation}\label{toda2}
\xymatrix{R&\S^{\abs{a}}R\ar[l]_{\bar{a}}&\S^{{\abs{a}}+{\abs{b}}}R\ar[l]|{\S^{\abs{a}} \bar{b}}^<(.85){\;}="c"\ar@/^30pt/[ll]^0_{\;}="d"&
\S^{{\abs{a}}+{\abs{b}}+{\abs{c}}}R\ar[l]^{\S^{{\abs{a}}+{\abs{b}}}\bar{c}}_<(.88){\;}="a"\ar@/_30pt/[ll]_0^{\;}="b"\ar@{=>}"a";"b"^f\ar@{=>}"c";"d"^e}.
\end{equation}
The pasting of this diagram is a self-track of the trivial map $0\colon \S^{{\abs{a}}+{\abs{b}}+{\abs{c}}}R\r R$.
Such a self-track is the same as a homotopy class
$$\S^{{\abs{a}}+{\abs{b}}+{\abs{c}}+1}R\To R,$$
which by definition is a generic element of the Toda bracket $\grupo{a,b,c}$.

The next definition encodes the secondary algebraic structure of the homotopy ring $\pi_*R$ endowed with the Toda
brackets. 

\begin{defn}\label{rmp}
Let $A$ be a ring and let $M$ be an $A$-bimodule. We say that $A$ has \emph{secondary oprations} with
coefficients in $M$ if there is given a bimodule homomorphism
$$\cdot\,\eta\colon A\otimes\Z/2\To M,$$
and for any three homogeneous elements $a,b,c\in A$ with $ab=0$ and $bc=0$ there is defined a coset (the
\emph{bracket operation})
$$\grupo{a,b,c}\subset M_{\abs{a}+\abs{b}+\abs{c}}\text{  of  }M_{\abs{a}+\abs{b}}\cdot c+a\cdot M_{\abs{b}+\abs{c}}$$
satisfying the following relations (whenever the brackets are defined):
\begin{enumerate}\renewcommand{\theenumi}{T\arabic{enumi}}
\item $0\in\grupo{a,b,c}$ provided $a$, $b$ or $c$ is zero.

\item $\grupo{a,b,c}$ is linear in each variable, i.e. 
\begin{eqnarray*}
\grupo{a+a',b,c}&\subset& \grupo{a,b,c}+\grupo{a',b,c},\\
\grupo{a,b+b',c}&=& \grupo{a,b,c}+\grupo{a,b',c},\\
\grupo{a,b,c+c'}&\subset& \grupo{a,b,c}+\grupo{a,b,c'}.
\end{eqnarray*}

\item $a\cdot\grupo{b,c,d}\subset \grupo{a\cdot b,c,d}$ and $\grupo{a,b,c}\cdot d\subset \grupo{a,b,c\cdot d}$.

\item $\grupo{a\cdot b,c,d}\subset\grupo{a,b\cdot c,d}\supset\grupo{a,b,c\cdot d}$.

\item $0\in\grupo{a,b,c}\cdot d+a\cdot\grupo{b,c,d}$,

\item $a\cdot \eta\in\grupo{2,a,2}$.
\end{enumerate}
\end{defn}

The \emph{desuspension} $\S^{-1}A$ of a ring $A$ is the $A$-bimodule with $(\S^{-1}A)_n=A_{1+n}$ and bimodule
structure
\begin{eqnarray*}
a\cdot (\S^{-1}b)\cdot c&=&(-1)^{\abs{a}}\S^{-1}(a\cdot b\cdot c).
\end{eqnarray*}
Here $a,b,c\in A$ are homogeneous elements and given $x\in A_n$ with $n\geq 1$ we denote by $\S^{-1}x$ to the
corresponding element in $(\S^{-1}A)_{n-1}$.
One can similarly define the desuspension of a right $A$-module, for which there are no signs involved in the
action.

\begin{thm}\label{main1}
Let $R$ be a connective ring spectrum. The ring $\pi_*R$ has secondary operations with coefficients in $\S^{-1}\pi_*R$, in the
sense of Definition \ref{rmp}. The homomorphism $\cdot\,\eta$ is defined by multiplication from the right
by the image of the stable Hopf map
$0\neq\eta\in\pi_1S\cong\Z/2$ under the ring homomorphism $\pi_*S\r\pi_*R$ induced by the unit $S\r R$ of the ring
spectrum, and the bracket operation is given by Toda brackets.
\end{thm}

This theorem follows from Theorems \ref{a2} and \ref{moin1} below.

\begin{rem}\label{coin1}
Alexander considered in \cite[Definition 2.1]{cmp} a notion of a ring with secondary operations similar to Definition \ref{rmp}.
Our relations (T1)--(T5) correspond to relations (1)--(5) in \cite[Definition 2.1]{cmp}
if $M=\S^{-1} A$. The homomorphism $\eta$, and therefore (T6) above, are not considered in \cite[Definition
2.1]{cmp}, although they appear in particular examples, see \cite[Theorem 6.2]{cmp}.
Relation (6) in \cite[Definition 2.1]{cmp} is not codified by
Definition \ref{rmp} since it is not a secondary relation, it has higher order. 
Alexander's relations are claimed in \cite{cmp} for some cobordism rings. These rings may arise as the homotopy groups of connective spectra, see
\cite[Section 4]{cmp}, which may be given the structure of ring spectra as in \cite[Example 4.15]{mtadsmc}.

Alexander's relations (1)--(6) coincide with the relations (3.5)--(3.8) in \cite{toda} previously proved by Toda
for the sphere spectrum.
\end{rem}

The homotopy groups $\pi_*M$ of a connective right $R$-module $M$ form a right $\pi_*R$-module which is also endowed
with Toda brackets 
$$\grupo{a,b,c}\in \pi_{\abs{a}+\abs{b}+\abs{c}+1}M,$$
defined for homogeneous elements $a\in \pi_*M$ and $b,c\in\pi_*R$ whenever $ab=0$ and $bc=0$, which is a coset of 
$$(\pi_{\abs{a}+\abs{b}+1}M)\cdot c+a\cdot (\pi_{\abs{b}+\abs{c}+1}R)$$
These Toda brackets are defined replacing
$R$  by $M$ on the left hand side of diagram (\ref{toda2}).

The folowing definition codifies the algebraic structure of Toda brackets in $\pi_*M$.

\begin{defn}\label{mmp}
Let $A$ be a ring with secondary operations with coefficients in the $A$-bimodule $M$ and let $N$ and $L$ be
right $A$-modules. We say that $N$ has \emph{secondary operations} with
coefficients $L$ if a right $A$-module homomorphism
$$\cdot\;\colon N\otimes_A M\To L$$
is given and
for any three homogeneous elements $a\in N$, $b,c\in A$ with $ab=0$ and $bc=0$ there is defined a coset (the
\emph{bracket operation})
$$\grupo{a,b,c}\subset L_{\abs{a}+\abs{b}+\abs{c}}\text{  of  }L_{\abs{a}+\abs{b}}c+aM_{\abs{b}+\abs{c}}$$
satisfying relations (T1)--(T5) in Definition \ref{rmp} for $a,a'\in N$ and $b,b',c,c',d\in A$.
\end{defn}

\begin{thm}\label{main2}
Let $R$ be a connective ring spectrum and let $K$ be a connective right $R$-module. 
The right $\pi_*R$-module $\pi_*K$ has secondary operations with coefficients in $\S^{-1}\pi_*K$, in the
sense of Definition \ref{mmp}. The homomorphism 
$$\cdot\;\colon \pi_*K\otimes_{\pi_*R} \S^{-1}\pi_*R\To \S^{-1}\pi_*K$$
is defined by the right $\pi_*R$-module structure of $\pi_*K$ according to the formula
\begin{eqnarray*}
m\cdot (\S^{-1}a)&=&(-1)^{\abs{m}}\S^{-1}(m\cdot a),
\end{eqnarray*}
and the bracket operation is given by Toda brackets.
\end{thm}

This theorem follows from Theorems \ref{m2} and \ref{moin2} below.

The homotopy groups of a commutative ring spectrum $\pi_*R$ form a commutative ring (in the graded sense) which carries,
appart from Toda brackets,
an additional operation called \emph{cup-one square}, 
$$\sq_1\colon\pi_{2n}R\To\pi_{4n+1}R,$$
defined as follows. 
Let $LR$ be a fibrant replacement of $R$ in the category of
all ring spectra. The ring spectrum $LR$ is no longer commutative, but it remains commutative up to a coherent 
track $\alpha_1$ (i.e. a homotopy class of homotopies) satisfying the idempotence and the hexagon axioms for symmetric monoidal
categories, compare \cite[Lemma 16.2]{2hg4}. Given $a\in\pi_{2n}R$ we take a representative $\bar{a}\colon
S^{2n}\r LR$ where the spectrum $S^m$ is the $m$-fold suspension of the sphere spectrum $S$, $S^m=\S^mS$. 
The symmetry isomorphism for the smash square of
an even-dimensional sphere $\tau_\wedge\colon S^{2n}\wedge S^{2n}\cong S^{2n}\wedge S^{2n}$ is homotopic to the
identity. 
We can choose a track $\hat{\tau}_{2n,2n}\colon \tau_{\wedge}\rr 1_{S^{2n}\wedge S^{2n}}$, there are two such choices. 
Consider the
following diagram where $\mu$ is the product in $LR$.
\begin{equation}\label{tcup1}
\xymatrix{S^{2n}\wedge S^{2n}\ar[dd]_{\;}="b"^{\tau_\wedge}\ar@/_30pt/[dd]_1^{\;}="a"
\ar[r]^-{\bar{a}\wedge\bar{a}}&LR\wedge
LR\ar[dd]_{\tau_\wedge}^{\;}="c"\ar[rd]_{\;}="d"^\mu&\\
&&LR\\
S^{2n}\wedge S^{2n}\ar[r]_-{\bar{a}\wedge\bar{a}}&LR\wedge LR\ar[ru]_\mu&
\ar@{<=}"a";"b"^{\hat{\tau}_{2n,2n}}\ar@{=>}"c";"d"_{\alpha_1}}
\end{equation}
The pasting of this diagram is a self-track of $\mu(\bar{a}\wedge\bar{a})$. The classical Barcus-Barratt-Rutter 
isomorphism allows us to identify this self-track 
with a homotopy class
$$Sq_1(a)\colon S^{4n+1}=\S(S^{2n}\wedge S^{2n})\To R$$
measuring the difference between the pasting of (\ref{tcup1}) and the identity self-track on $\mu(\bar{a}\wedge\bar{a})$.
This element $Sq_1(a)\in\pi_{4n+1}R$ is the cup-one square of $a$. One can check that $Sq_1(a)$ does not depend
on the representative $\bar{a}$. However in general it does depend on the choice of $\hat{\tau}_{2n,2n}$. The difference
bewteen the two possible definitions of $\sq_1$, depending on the choice of $\hat{\tau}_{2n,2n}$, is computed in
\cite[Lemma 9.11]{2hg4}, see Lemma \ref{lor2}.

The relations between Toda brackets and cup-one squares in the homotopy groups of a commutative ring spectrum is
algebraically encoded by the following definition.

\begin{defn}\label{crmp}
A commutative ring $A$ with \emph{commutative secondary operations} with coefficients in an $A$-module $M$ is a ring
with secondary operations in the sense of Definition \ref{rmp} together with maps
$$\sq_1\colon A_{2n}\To M_{4n},\;\;n\geq 0,$$
such that the following further axioms hold:
\begin{enumerate}\renewcommand{\theenumi}{T\arabic{enumi}}\setcounter{enumi}{6}
\item $\grupo{a,b,c}=(-1)^{\abs{a}\abs{b}+\abs{b}\abs{c}+\abs{c}\abs{a}+1}\grupo{c,b,a}$.

\item $0\in (-1)^{\abs{a}\abs{c}}\grupo{a,b,c}+(-1)^{\abs{b}\abs{a}}\grupo{b,c,a}
+(-1)^{\abs{c}\abs{b}}\grupo{c,a,b}$,

\item for $\abs{a}$ odd $\grupo{a,b,a}\cap(-1)^{\abs{a}\abs{b}}\grupo{b,a,2a}\neq\emptyset$,

\item for $\abs{a}$ even $(-1)^{\abs{a}\abs{b}}b\cdot\sq_1(a)\in\grupo{a,b,a}$,

%\item $2\cdot\sq_1(a)=\frac{\abs{a}}{2}a^2\cdot\eta$,

\item $\sq_1(a+b)=\sq_1(a)+\sq_1(b)+\left(\frac{\abs{a}}{2}+1\right)\cdot a\cdot b\cdot\eta$,

\item $\sq_1(a\cdot b)=a^2\cdot\sq_1(b)+\sq_1(a)\cdot b^2+\frac{\abs{a}\abs{b}}{4}\cdot a^2\cdot
b^2\cdot\eta$.
\end{enumerate}
\end{defn}

\begin{thm}\label{main3}
Let $R$ be a connective commutative ring spectrum. The ring $\pi_*R$ has commutative secondary operations with
coefficients in $\S^{-1}\pi_*R$, in the
sense of Definition \ref{crmp}. The operation $\sq_1$ is the cup-one square for an explicit choice of tracks
$\hat{\tau}_{2n,2n}$, $n\geq 0$, and the rest of the structure is
given by Theorem \ref{main1}.
\end{thm}

This theorem follows from Theorems \ref{ca2} and \ref{moin3} below.

\begin{rem}\label{coin2}
There is also a notion of a commutative ring with commutative secondary operations in \cite[Definition 2.1]{cmp}.
This notion however does not codify the operation $\sq_1$.
Our relations (T7) and (T8) correspond to relations (7) and (8) in \cite[Definition 2.1]{cmp}
if $M=\S^{-1} A$. These relations are claimed in \cite{cmp} for some commutative cobordism rings which may
arise as the homotopy groups of commutative connective ring spectra, see Remark \ref{coin1}. 
Notice that there is a misprint in the exponent of $(-1)$ in relation (8) of \cite[Definition 2.1]{cmp}. It
does not include the summand $+1$. This misprint does not appear in Toda's relations for the case of the sphere spectrum, see (3.9)
in \cite{toda}. Relation (T9) corresponds to the first half of Toda's (3.10) in \cite{toda}. The second half is a
weak version of (T10) which avoids the use of $\sq_1$. 

$H_\infty$-ring spectra in the sense of \cite{hirs} are an early version of commutative ring spectra ``up to
homotopy". The operations $\sq_1$ are closely related to the power operations for $H_\infty$-ring spectra
considered in \cite[V.1]{hirs}. More precisely, the operation $P^{n+1}$ in \cite[V]{hirs} on $\pi_n$ for $p=2$ and
$n=2k$ corresponds to the set $P^{n+1}(a)=\set{\sq_1(a),\sq_1(a)+a^2\cdot\eta}$. Then relation (T11) above
implies the deviation from additivity indicated in \cite[V Table 1.3]{hirs} and relation (T12) implies the first
equation of \cite[Proposition V.1.10]{hirs} and gives the explicit value for the constant $c_{n,m}$ which is not
determined in \cite{hirs}. Maybe one of the most surprising implications of Theorem \ref{main3} is the existence of
choices $\sq_1(a)\in P^{n+1}(a)$, for $p=2$, $n=2k$, and $k\geq0$, satisfying relations (T11) and (T12).

Any commutative ring $A$ with commutative secondary operations with coefficients in an $A$-module $M$ has
the following remarkable property. If we define $\sq^\omega_1(a)=\sq_1(a)+a^2\cdot\eta$, then
$\sq^\omega_1$
also satisifies the axioms in Definition \ref{crmp}.
\end{rem}

In the following proposition we record some additional relations between cup-one squares derived from Definition \ref{crmp}.

\begin{prop}\label{mas}
Let $A$ be a commutative ring with commutative secondary operations with coefficients in the $A$-module $M$. Then
\begin{enumerate}
\item $\sq_1(1)=0$,

\item $\sq_1(2)=1\cdot\eta$,

\item $2\cdot\sq_1(a)=\frac{\abs{a}}{2}\cdot a^2\cdot\eta$,

\item $\sq_1(2\cdot a)=a^2\cdot\eta$.
\end{enumerate}
\end{prop}

\begin{proof}
Equation (1) follows from (T12) applied to $a=b=1$, and (2) follows from (T11) and (1). Applying (T11) to $a+a$ and
(T12) to $2\cdot a$ we obtain the equation
\begin{eqnarray*}
2\cdot\sq_1(a)+\left(\frac{\abs{a}}{2}+1\right)\cdot a^2\cdot\eta&=&4\cdot\sq_1(a)+a^2\cdot\eta.
\end{eqnarray*}
Here we use (2) to identify $\sq_1(2)\cdot a^2=a^2\cdot\eta$. Equation (3) follows from this one. Finally (4)
follows from (T11) and (3).
\end{proof}

Similar relations are shown in \cite[V.1]{hirs} for the power operations on the homotopy groups of
$H_\infty$-ring spectra.

\begin{rem}
By Proposition \ref{mas} (2) the structure homomorphism $\cdot\,\eta$ of a commutative ring with commutative
secondary operations is determined by the operation $\sq_1$, so one could restate Definition \ref{crmp} just in
terms of the bracket and $\sq_1$.
\end{rem}

\section{Quadratic pair modules}
\renewcommand{\theequation}{M\arabic{equation}}\setcounter{equation}{0}

The topological theorems of this paper are proved by using the quadratic algebraic models for ring and module
spectra defined in \cite{2hg4}. In this section we recall the basics on the necessary quadratic algebra, see
\cite{ecg, qaI}.

A \emph{quadratic pair module} $C$ is a diagram 
\begin{equation*}
\xymatrix{&C_{ee}\ar[ld]_P&\\C_{1}\ar[rr]_\partial&&C_{0}\ar[lu]_H}
\end{equation*}
where $C_0$ and $C_1$ are groups, $C_{ee}$ is an abelian group, $P$ and $\partial$ are homomorphisms, and $H$ is
a \emph{quadratic map}, i.e. the \emph{crossed efect}
$$(x_1|x_2)_H=H(x_1+x_2)-H(x_2)-H(x_1),\;\;x_i\in C_0,$$
is bilinear. Moreover, the following equations hold for $x,x_i\in C_0$, $s_i\in C_1$ and $a\in C_{ee}$.
\begin{eqnarray}
\label{qpm1} PH\partial P(a)&=&P(a)+P(a),\\
\label{qpm2} H(x+\partial P(a))&=&H(x)+H\partial P(a),\\
\label{qpm3} PH(\partial(s_1)+\partial(s_2))&=&PH\partial(s_1)+PH\partial(s_2)+[s_1,s_2],\\
\label{qpm4} \partial PH(x_1+x_2)&=&\partial PH(x_1)+\partial PH(x_2)+[x_1,x_2],
\end{eqnarray}
see \cite[2.4]{3mlc}.
Here $[\alpha,\beta]=-\alpha-\beta+\alpha+\beta$ denotes the commutator bracket of two elements $\alpha,\beta \in
G$ in a group $G$.

It follows from the axioms that $C_0$ and $C_1$ are groups of nilpotency class $2$, so commutators are central
and bilinear, $\partial(C_1)$ is a normal
subgroup of $C_0$, and $P$ and $\ker\partial$ are central. The quadratic map $H$ satisfies
\begin{eqnarray*}
H(0)&=&0,\\
H(-x)&=&-H(x)+(x|x)_H.
\end{eqnarray*}
For any quadratic pair module the function 
$$T=H\partial P-1\colon X_{ee}\To X_{ee}$$ is an involution, i.e. a homomorphism with $T^2=1$. Using (\ref{qpm4})
cone can check that
\begin{eqnarray*}
T(x_1|x_2)_H=-(x_2|x_1)_H.
\end{eqnarray*}
By (\ref{qpm1}) $T$ satisfies
$PT=P$, therefore
\begin{eqnarray*}
P(x_1|x_2)_H=-P(x_2|x_1)_H.
\end{eqnarray*}
Moreover, $$\Delta\colon C_0\To C_{ee}\colon x\mapsto (x|x)_H-H(x)+TH(x)$$
is a homomorphism which satisfies %$T\Delta=-\Delta$ and 
$P\Delta(x)=P(x|x)_H$.

The \emph{homology} of a quadratic pair module $C$ is given by the abelian groups defined as
\begin{eqnarray*}
h_0 C&=&\coker\left(\partial\colon C_1\r C_0\right),\\
h_1 C&=&\ker\left(\partial\colon C_1\r C_0\right).
\end{eqnarray*}
The \emph{$k$-invariant} of $C$ is the natural homomorphism
\begin{equation*}
\cdot\;\eta\colon h_0C\otimes\Z/2\To h_1C
\end{equation*}
given by the formula 
\begin{eqnarray*}
x\cdot\eta &=&P(x|x)_H\quad = \quad P\Delta(x).
\end{eqnarray*}

A \emph{morphism} $f\colon C\r D$ of quadratic pair modules is given by three homomorphisms $f_i\colon C_i\r D_i$,
$i=0,1,ee$, commuting with the structure maps, i.e. $f_0\partial=\partial f_1$, $f_1P=Pf_{ee}$, $f_{ee}H=Hf_0$.
Morphisms of quadratic pair modules are also compatible with $T$, $\Delta$, and $\cdot\,\eta$. A
\emph{quasi-isomorphism} is a morphism inducing isomorphisms in $h_0$ and $h_1$.

\section{Quadratic pair algebras}
\renewcommand{\theequation}{A\arabic{equation}}\setcounter{equation}{0}

In this section we recall the nature of the quadratic algebraic models for ring spectra constructed in
\cite{2hg4} and we prove Theorem \ref{main1}.

A \emph{quadratic pair algebra} is an $\N$-graded quadratic pair module 
$B=\{B_{n,*}\}_{n\in \N}$, together with multiplications, $n,m\in \N$,
$$\begin{array}{c}
B_{n,0}\times B_{m,0}\st{\cdot}\To B_{n+m,0},\\
B_{n,0}\times B_{m,1}\st{\cdot}\To B_{n+m,1},\\
B_{n,1}\times B_{m,0}\st{\cdot}\To B_{n+m,1},\\
B_{n,ee}\times B_{m,ee}\st{\cdot}\To B_{n+m,ee},
\end{array}$$
and an element $1\in B_{0,0}$ with $H(1)=0$ which is a (two-sided) unit for the
first three multiplications and such that $(1|1)_H\in B_{0,ee}$ is a (two-sided) unit for the fourth
multiplication. These multiplications are associative in all possible ways.
Moreover, the following lists
of equations are satisfied for $x,x_i\in B_{*,0}$, $s,s_i\in
B_{*,1}$ and $a_i\in B_{*,ee}$.
The multiplications $\cdot$ are
always right linear
\begin{eqnarray}
\label{rl} x_1\cdot(x_2+x_3)&=&x_1\cdot x_2 + x_1\cdot x_3,\\
\nonumber x\cdot(s_1+s_2)&=&x\cdot s_1 + x\cdot s_2,\\
\nonumber s\cdot(x_1+x_2)&=&s\cdot x_1 + s\cdot x_2,\\
\nonumber a_1\cdot(a_2+a_3)&=&a_1\cdot a_2 + a_1\cdot a_3.
\end{eqnarray}
The multiplications $\cdot$ satisfy the following left
distributivity laws
\begin{eqnarray}
\label{ldl} (x_1+x_2)\cdot x_3&=&x_1\cdot x_3+x_2\cdot x_3 + \partial P((x_2|x_1)_H\cdot H(x_3)),\\
\nonumber (x_1+x_2)\cdot s&=&x_1\cdot s+x_2\cdot s + P((x_2|x_1)_H\cdot H\partial(s)),\\
\nonumber (s_1+s_2)\cdot x&=&s_1\cdot x+s_2\cdot x + P((\partial(s_2)|\partial(s_1))_H\cdot H(x)),\\
\nonumber (a_1+a_2)\cdot a_3&=&a_1\cdot a_3+a_2\cdot a_3.
\end{eqnarray}
The homomorphisms $\partial$ are compatible with the multiplications
$\cdot$ in the following sense
\begin{eqnarray}
\label{dm} \partial(x\cdot s)&=&x\cdot\partial(s),\\
\nonumber \partial(s\cdot x)&=&\partial(s)\cdot x,\\
\nonumber \partial(s_1)\cdot s_2&=&s_1\cdot\partial(s_2).
\nonumber \end{eqnarray}
And finally, we have compatibility conditions for the
multiplications $\cdot$ and the maps $P$, $H$, $\Delta$, and
$(-|-)_H$,
\begin{eqnarray}
\label{pml} P((x|x)_H\cdot a)&=&x\cdot P(a),\\
\label{pmr} P(a\cdot\Delta(x))&=&P(a)\cdot x,\\
\label{hm} H(x_1\cdot x_2)&=& (x_1|x_1)_H\cdot H(x_2)+H(x_1)\cdot\Delta(x_2),\\
\label{hdpm} H\partial P(a_1\cdot a_2)&=&H\partial P(a_1)\cdot a_2+a_1\cdot H\partial P(a_2)\\
\nonumber &&-H\partial P(a_1)\cdot
H\partial P(a_2),\\
\label{cem} (x_1\cdot x_2|x_3\cdot x_4)_H&=&(x_1|x_3)_H\cdot (x_2|x_4)_H.
\end{eqnarray}

Ungraded quadratic pair algebras were first considered in \cite{3mlc} in order to represent classes in third
Mac Lane cohomology. The graded notion, which is the one we mainly use in this paper, was introduced in
\cite{2hg4}.

A \emph{morphism} of quadratic pair algebras is a morphism of graded quadratic pair modules which preserves the
products $\,\cdot\,$. A \emph{quasi-isomorphism} is a morphism inducing isomorphisms in $h_0$ and $h_1$.

\setcounter{equation}{0}
\numberwithin{equation}{section}
\renewcommand{\theequation}{\thesection.\arabic{equation}}

For $B$ a quadratic pair algebra $h_0B$ is a ring ($\N$-graded) and $h_1B$ is an $h_0B$-bimodule in a natural way.
Moreover, the $k$-nvariant 
\begin{equation}\label{qpak}
\cdot\;\eta\colon h_0B\To h_1B
\end{equation}
is an $h_0B$-bimodule homomorphism by (\ref{pml},\ref{pmr},\ref{cem}).

The relations in the following lemma are consequences of (\ref{ldl}).

\begin{lem}\label{menos}
With the notation above the following equations hold.
\begin{enumerate}
\item $0\cdot x_2=0$,
\item $(-x_1)\cdot x_2=-x_1\cdot x_2+\partial P((x_1|x_1)_H\cdot H(x_2))$,
\item $(-x)\cdot s=-x\cdot s+P((x|x)_H\cdot H\partial (s))$.
\end{enumerate}
\end{lem}

\begin{defn}\label{qpamas}
Let $B$ be a quadratic pair algebra. Given elements $a,b,c\in h_0 B$,
of degree $p,q,r\in\N$ with $ab=0$ and $bc=0$ the
\emph{Massey product}  is the subset
$$\grupo{a,b,c}\subset h_1B_{p+q+r},$$
which is a coset of the subgroup $$(h_1B_{p+q})\cdot c+a\cdot(h_1B_{q+r}),$$
defined as follows. Given $\bar{a}\in B_{p,0}$, $\bar{b}\in
B_{q,0}$, $\bar{c}\in B_{r,0}$ representing $a$, $b$, $c$, there exist $\overline{ab}\in B_{p+q,1}$,
$\overline{bc}\in B_{q+r,1}$ such
that $\partial(\overline{ab})=\bar{a}\cdot\bar{b}$,
$\partial(\overline{bc})=\bar{b}\cdot\bar{c}$ and one can easily check that 
$$-\overline{ab}\cdot\bar{c}+\bar{a}\cdot \overline{bc}\in h_1B_{p+q+r}\subset B_{p+q+r,1}.$$
The coset $\grupo{a,b,c}\subset h_1B_{p+q+r}$ 
coincides with the set of elements obtained in this way for all
different choices of $\bar{a}$, $\bar{b}$, $\bar{c}$, $\overline{ab}$ and $\overline{bc}$.
\end{defn}

In \cite{2hg4} we prove the following theorem as a main result.

\begin{thm}[{\cite[Theorem 6.4]{2hg4}}]\label{a2}
There is a functor
$$\pi_{*,*}\colon \left(
\text{\emph{connective ring spectra}}\right)\To\left(\text{\emph{quadratic pair algebras}}\right)$$
together with natural isomorphisms
\begin{eqnarray*}
h_0\pi_{*,*} R&\cong&\pi_*R, \text{ of rings},\\
h_1\pi_{*,*} R&\cong&\S^{-1}\pi_*R, \text{ of bimodules},
\end{eqnarray*}
such that the Massey products in $\pi_{*,*} R$ coincide with the Toda brackets in $\pi_*R$. Moreover, using the
isomorphisms as identifications the algebraically-defined
$k$-invariant of the quadratic pair algebra $\pi_{*,*} R$ 
$$\cdot\;\eta\colon \pi_*R\otimes\Z/2\To\S^{-1}\pi_*R,$$
coincides with the multiplication by the image of the stable Hopf map under the homomorphism $\pi_*S\r \pi_*R$
induced by the unit $S\r R$.
\end{thm}

Theorem \ref{main1} will then follow from the following one.

\begin{thm}\label{moin1}
If $B$ is a quadratic pair algebra then the $k$-invariant (\ref{qpak}) and 
the Massey products in Definition \ref{qpamas} endow $h_0B$ with the structure of a
ring with secondary operations with coefficients in $h_1B$ in the sense of Definition \ref{rmp}.
\end{thm}

For the sake of simplicity in the proof of Theorem \ref{moin1} we will use assume that $B$ satisfies the property (H).

\bigskip

\begin{itemize}
\item[(H)] Any element in $x\in h_0 B$
is the image of an element $\bar{x}\in {B}_{*,0}$ with $H(\bar{x})=0$.
\end{itemize}

\bigskip

This property is not unusual. For instance, given a ring spectrum $R$ the quadratic pair algebra $\pi_{*,*}R$
defined by Theorem \ref{a2} satisfies property (H). Indeed the following lemma holds.

\begin{lem}\label{kerh}
Given a quadratic pair algebra $B$ there is another one $\widehat{B}$ satisfying property (H) and a natural quasi-isomorphism $B\r
\widehat{B}$.
\end{lem}

Here $\widehat{B}$ is a fibrant replacement of $B$ in the cofibration category of quadratic pair algebras and is
obtained ``attaching cells" to $B$. We will not discuss the homotopical aspects of quadratic pair algebras in
this paper, so we leave the proof of Lemma \ref{kerh} to the interested reader. This lemma shows that there is no
loss of generality if we only prove Theorem \ref{moin1} for quadratic pair algebras satisfying property (H).

\begin{rem}\label{nancy}
Before beginning the proof of Theorem \ref{moin1} we want to remark that in orther to check the inclusions and
equalities in Definition \ref{rmp} it is enough to check that the brackets have an alement in common. Then the
inclusion (resp. equality) follows from the obvious analogous inclusion (resp. equality) between the
indeterminacies which is clear in all cases.
The same applies to the proof of Theorem \ref{moin3}.
\end{rem}

\begin{proof}[Proof of Theorem \ref{moin1}]
We assume that all representatives chosen in $B_{*,0}$ are in $\ker H$. Let us check that equations (T1)--(T6)
hold.

(T1) If $a=0$ we can take $\bar{a}=0$ and $\overline{ab}=0$ so $-\overline{ab}\cdot\bar{c}+\bar{a}\cdot
\overline{bc}=0$. Similarly in the other two cases.

(T2) We can take $\overline{a+a'}=\bar{a}+\bar{a}'$ and by (\ref{ldl}) we can also take
$\overline{(a+a')b}=\overline{ab}+\overline{a'b}$, therefore
\begin{eqnarray*}
-\overline{(a+a')b}\cdot\bar{c}+\overline{a+a'}\cdot\overline{bc}\!\!\!\!\!\!&\st{\mbox{\scriptsize
(\ref{ldl})}}=&\!\!\!\!\!\!
-\overline{a'b}\cdot\bar{c}-\overline{ab}\cdot\bar{c}+\bar{a}\cdot\overline{bc}+\bar{a}'\cdot\overline{bc}-P((a'|a)_H\cdot
\underbrace{H(\bar{b}\cdot\bar{c})}_{\mbox{\scriptsize (\ref{hm})}\;\;=0})\\
\mbox{\scriptsize
(\ref{qpm3})}\quad&=&\!\!\!\!\!\!-\overline{a'b}\cdot\bar{c}+\bar{a}\cdot\overline{bc}-\overline{ab}\cdot\bar{c}+\bar{a}'\cdot\overline{bc}
-\underbrace{P(\bar{a}\cdot\bar{b}\cdot\bar{c}|\bar{a}\cdot\bar{b}\cdot\bar{c})_H}_{=\;a\cdot b\cdot c\cdot
\eta\;=\;0}.
\end{eqnarray*}
One proceeds similarly with the two other variables.

(T3) By (\ref{dm}) $\partial(\bar{a}\cdot\overline{bc})=\bar{a}\cdot\bar{b}\cdot\bar{c}$, hence the first
equation in (T3) follows from 
\begin{eqnarray*}
\bar{a}\cdot(-\overline{bc}\cdot\bar{d}+\bar{b}\cdot\overline{cd})&\st{\mbox{\scriptsize (\ref{rl})}}=&
-(\bar{a}\cdot\overline{bc})\cdot\bar{d}+(\bar{a}\cdot\bar{b})\cdot\overline{cd}.
\end{eqnarray*}
The second one follows similarly.

(T4) By (\ref{dm}) $\partial(\bar{b}\cdot\overline{cd})=\bar{b}\cdot\bar{c}\cdot\bar{d}$ therefore
$-\overline{abc}\cdot\bar{d}+\bar{a}\cdot\bar{b}\cdot\overline{cd}$ lies in both $\grupo{ab,c,d}$ and $\grupo{a,bc,d}$. Similarly the element $-\overline{ac}\cdot\bar{c}\cdot\bar{d}
+\bar{a}\cdot\overline{bcd}$ belongs to the other two Massey products.

(T5) This follows from
\begin{eqnarray*}
(-\overline{ab}\cdot\bar{c}+\bar{a}\cdot\overline{bc})\cdot\bar{d}
+\bar{a}\cdot(-\overline{bc}\cdot\bar{d}+\bar{b}\cdot\overline{cd})&\st{\mbox{\scriptsize (\ref{rl},\ref{ldl})}}=&
-\overline{ab}\cdot\bar{c}\cdot\bar{d}+\bar{a}\cdot\overline{bc}\cdot\bar{d}\\
&&-\bar{a}\cdot\overline{bc}\cdot\bar{d}+\bar{a}\cdot\bar{b}\cdot\overline{cd}\\
&=&-\overline{ab}\cdot\bar{c}\cdot\bar{d}+\bar{a}\cdot\bar{b}\cdot\overline{cd}\\
&=&-\overline{ab}\cdot\partial(\overline{cd})+\partial(\overline{ab})\cdot\overline{cd}\\
\mbox{\scriptsize (\ref{dm})}&=&0.
\end{eqnarray*}

Finally (T6) follows from \cite[Proposition 6.6]{2hg4}.
\end{proof}

\section{Modules over quadratic pair algebras}

The quadratic algebraic models of module spectra leading to Theorem \ref{main2} are as follows.

Let $B$ be a quadratic pair algebra. A \emph{right $B$-module} 
is an $\N$-graded quadratic pair
module $M=\{M_{n,*}\}_{n\in \N}$ together with multiplications, $n,m\geq 0$,
$$\begin{array}{c}
M_{n,0}\times B_{m,0}\st{\cdot}\To M_{n+m,0},\\
M_{n,0}\times B_{m,1}\st{\cdot}\To M_{n+m,1},\\
M_{n,1}\times B_{m,0}\st{\cdot}\To M_{n+m,1},\\
M_{n,ee}\times B_{m,ee}\st{\cdot}\To M_{n+m,ee}.
\end{array}$$
These multiplications are associative with respect to the multiplications in $B$.
Moreover, $1\in B_{0,0}$ acts trivially on $M_{*,0}$ and $M_{*,1}$, and  $(1|1)_H\in B_{0,ee}$ acts trivially on
$M_{*,ee}$. Furthermore, equations (\ref{rl})--(\ref{cem}) hold when we replace the elements on the left of
any multiplication $\cdot$ by elements in $M$.

If $M$ is a right $B$-module then $h_0M$ and $h_1M$ are right $h_0B$-modules and there is a natural right $h_0B$-module
homomorphism
\begin{equation}\label{pairing}
h_0M\otimes_{h_0B}h_1B\st{\cdot}\To h_1M,
\end{equation}
see \cite[7]{2hg4},
extending the $k$-invariant since $x\cdot\eta=x\cdot P(1|1)_H$ for $x\in h_0M$ by (\ref{pml}). The $k$-invariant
of a right $B$-module $M$ is a right $h_0B$-module homomorphism by (\ref{pml},\ref{pmr},\ref{cem}).

\begin{defn}\label{modmas}
Given a right $B$-module  $M$ and elements $a\in h_0M$, $b,c\in h_0B$,
of degree $p,q,r,$ such that $ab=0$ and $bc=0$ there is defined a \emph{Massey product}
$$\grupo{a,b,c}\subset h_1M_{p+q+r}$$
by the same procedure as in Definition \ref{qpamas}
which is a coset of
$$(h_1M_{p+q})\cdot c+a\cdot(h_1B_{q+r}).$$
\end{defn}

We show the following theorem in \cite{2hg4}.

\begin{thm}[{\cite[Theorem 7.4]{2hg4}}]\label{m2}
Let $R$ be a connective ring spectrum. There is a functor
$$\pi_{*,*}\colon\left(\text{\emph{connective right $R$-modules}}\right)\To
\left(\text{\emph{right $(\pi_{*,*}R)$-modules}}\right).$$
Here the quadratic pair algebra $\pi_{*,*}R$ is obtained by using the functor in Theorem \ref{a2}.
Moreover, if we use the isomorphisms in Theorem \ref{a2} as identifications
then for any right $R$-module $K$ there are natural isomorphisms of right $\pi_*R$-modules
\begin{eqnarray*}
h_0\pi_{*,*}K&\cong&\pi_*K,\\
h_1\pi_{*,*}K&\cong&\S^{-1}\pi_*K.
\end{eqnarray*} 
Using these
isomorphisms as identifications the algebraically-defined
homomorphism (\ref{pairing}) associated to $\pi_{*,*} K$
$$\cdot\;\colon \pi_*K\otimes_{\pi_*R} \S^{-1}\pi_*R\To \S^{-1}\pi_*K$$
is defined by the right right $\pi_*R$-module structure of $\pi_*K$ according to the formula
\begin{eqnarray*}
m\cdot (\S^{-1}a)&=&(-1)^{\abs{m}}\S^{-1}(m\cdot a).
\end{eqnarray*}
In particular the $k$-invariant of $\pi_{*,*}K$ coincides with the multiplication by the stable Hopf map $\eta$.
Furthermore, Massey products in $\pi_{*,*}K$
coincide with Toda brackets in $\pi_*K$.
\end{thm}

Now Theorem \ref{main2} follows from the following theorem.

\begin{thm}\label{moin2}
If $B$ is a quadratic pair algebra and $M$ is a right $B$-module then (\ref{pairing}) and the Massey products in Definition \ref{modmas}
endow $h_0M$ with the structure of a
module with secondary operations with coefficients in $h_1M$ in the sense of Definition \ref{mmp}.
\end{thm}

The proof is completely analogous to the proof of Theorem \ref{moin1} so we leave it to the reader.

\section{Symmetric track groups}
\renewcommand{\theequation}{S\arabic{equation}}\setcounter{equation}{0}

In order to describe the quadratic algebraic models associated to commutative ring spectra we need to endow
quadratic pair modules with symmetries coming from the action of sign groups as we describe in this section.

A \emph{sign group} is a diagram in the category of groups
\begin{equation*}
\set{\pm1}\st{\imath}\hookrightarrow G_\vc\st{\delta}\twoheadrightarrow G\st{\varepsilon}\To \set{\pm1}
\end{equation*}
where the first two morphisms form an extension. By abuse of notation we denote this sign group just by $G_\vc$.
The group law of the groups defining a sign group is denoted multiplicatively.

Given a sign group $G_\vc$ the \emph{``group ring''} $A(G_\vc)$ is the ungraded 
quadratic pair algebra with
generators
\begin{itemize}
\item $[g]$ for any $g\in G$ on the $0$-level,
\item $[t]$ for any $t\in G_\vc$ on the $1$-level,
\item no generators on the $ee$-level,
\end{itemize}
satisfying the following relations for $g,h\in G$, $s,t\in G_\vc$ and $\omega=\imath(-1)$.
\begin{eqnarray}
\label{s1} H[g]&=&0,\\{}
\label{s2} [1]&=&1\text{ for } 1\in G,\\{}
\label{s3} [gh]&=&[g]\cdot[h],\\{}
\label{s4} \partial[t]&=&-[\delta(t)]+\varepsilon\delta(t),\\{}
\label{s5} [st]&=&[\delta(s)]\cdot[t]+\varepsilon\delta(t)\cdot[s]
+\binom{\varepsilon\delta(s)}{2}\binom{\varepsilon\delta(t)}{2}P(1|1)_H,\\{}
\label{s6} [\omega]&=&P(1|1)_H.
\end{eqnarray}

Sign groups were introduced in \cite{2hg2}, and the ``group ring'' of a sign group was first considered in
\cite{2hg3}.

\renewcommand{\theequation}{\thesection.\arabic{equation}}\setcounter{equation}{0}

The following lemma follows easily from (\ref{s1},\ref{s3}) and the fact that $H$ is quadratic.

\begin{lem}\label{hdt}
Given $t\in G_\vc$ the following equation  holds.
\begin{eqnarray*}
H\partial[t]&=&([\delta(t)]|-\partial[t])_H+\binom{\varepsilon\delta(t)}{2}(1|1)_H.
\end{eqnarray*}
\end{lem}

The following useful relation follows from (\ref{s2},\ref{s5}).
\begin{lem}\label{10}
For $1\in G_\vc$ we have $[1]=0$.
\end{lem}

The action of $A_0(G_\vc)$ on the left of $A_1(G_\vc)$ is dertermined by relation (\ref{s5}) in terms of the
group structure of $G_\vc$ since $\delta$ is surjective. The right action is given by the following lemma.

\begin{lem}\label{ader}
For $s,t\in G_\vc$ the following relation holds in $A(G_\vc)$.
\begin{eqnarray*}
[st]\!\!\!\!\!&=&\!\!\!\!\![s]\cdot[\delta(t)]+\varepsilon\delta(s)\cdot[t].
%+\binom{\varepsilon\delta(s)}{2}\binom{\varepsilon\delta(t)}{2}P(1|1)_H.
\end{eqnarray*}
\end{lem}

\begin{proof}
 On one hand by (\ref{s4},\ref{rl},\ref{menos}.3)
\begin{eqnarray*}
[s]\cdot\partial[t]&=&-[s]\cdot[\delta(t)]+\varepsilon\delta(t)\cdot[s]-\binom{\varepsilon\delta(t)}{2}PH\partial[s].
\end{eqnarray*}
On the other hand by (\ref{s4},\ref{ldl},\ref{menos}.3,\ref{hdt},\ref{cem},\ref{qpm3},\ref{dm})
%\begin{eqnarray*}
%\partial[s]\cdot[t]&=&-[\delta(s)]\cdot[t]+\varepsilon\delta(s)\cdot[t]\\
%&&+P((\overbrace{(\varepsilon\delta(s)|-[\delta(s)])_H+([\delta(s)]|[\delta(s)])_H}^{=(-\partial[s]|[\delta(s)])_H}+\binom{\varepsilon\delta(s)}{2}(1|1)_H)\\
%&&\;\;\;\cdot(([\delta(t)]|-\partial[t])_H+\binom{\varepsilon\delta(t)}{2}(1|1)_H))\\
%\mbox{\scriptsize(\ref{cem},\ref{qpm3})}&=&-[s]\cdot[\delta(t)]-[\delta(s)]\cdot[t]+[s]\cdot[\delta(t)]+\varepsilon\delta(s)\cdot[1]
%\end{eqnarray*}
\begin{eqnarray*}
\partial[s]\cdot[t]&=&-[\delta(s)]\cdot[t]+\varepsilon\delta(s)\cdot[t]
+P((-\partial[s]|[\delta(s)])_H\cdot H\partial[t])\\
&=&-[s]\cdot[\delta(t)]-[\delta(s)]\cdot[t]+[s]\cdot[\delta(t)]+\varepsilon\delta(s)\cdot[t]\\
&&+\binom{\varepsilon\delta(t)}{2}P(-\partial[s]|[\delta(s)])_H.
\end{eqnarray*}
By (\ref{dm}) $\partial[s]\cdot[t]=[s]\cdot\partial[t]$, hence by using the two previous equations together with
(\ref{hdt}) one obtains
\begin{eqnarray*}
[\delta(s)]\cdot[t]+\varepsilon\delta(t)\cdot[s]+\binom{\varepsilon\delta(s)}{2}\binom{\varepsilon\delta(t)}{2}P(1|1)_H&=&[s]\cdot[\delta(t)]+\varepsilon\delta(s)\cdot[t].
\end{eqnarray*}
Now the lemma follows from (\ref{s5}).
\end{proof}

%\begin{lem}\label{mame}
%Given an $A(G_\vc)$-module $C$, $x\in\ker\partial\subset C_1$ and $g\in G$we have $x\cdot[g]=\varepsilon(g) x$.
%\end{lem}

%\begin{proof}
%Let $t\in G_\vc$ be an element with $\delta(t)=g$.
%\begin{eqnarray*}
%0&=&\partial(x)\cdot[t]\\
%\mbox{\scriptsize (\ref{dm})}\qquad&=&x\cdot\partial[t]\\
%\mbox{\scriptsize (\ref{s4})}\qquad&=&x\cdot(-[g]+\varepsilon(g))\\
%\mbox{\scriptsize (\ref{rl})}\qquad&=&-x\cdot[g]+\varepsilon(g)x.
%\end{eqnarray*}
%\end{proof}

The homology of ``group rings'' of sign groups can be easily computed.

\begin{lem}\label{mame}
There are natural isomorphisms 
\begin{eqnarray*}
h_0A(G_\vc)&\cong&\Z,\\
h_1A(G_\vc)&\cong&\Z/2.
\end{eqnarray*}
The first one is induced by $[g]\mapsto \varepsilon(g)$, and $h_1A(G_\vc)$ is generated by $[\omega]$.
\end{lem}

The main examples of sign groups are the \emph{symmetric track groups} $\symt{n}$ associated
to the sign homomorphism of the symmetric groups $$\varepsilon=\sign\colon\sym{n}\r\set{\pm1}.$$ The group $\symt{n}$ has a presentation with
generators $\omega$, $t_i$, $1\leq i\leq n-1$,
and relations
\begin{eqnarray}
\nonumber t_i^2&=&1\text{ for }1\leq i\leq n-1,\\
\label{elados} (t_it_{i+1})^3&=&1\text{ for }1\leq i\leq n-2,\\
\nonumber \omega^2&=&1,\\
\nonumber t_i\omega&=&\omega t_i\text{ for }1\leq i\leq n-1,\\
\nonumber t_it_j&=&\omega t_jt_i\text{ for }1\leq i<j-1\leq n-1.
\end{eqnarray}
Moreover, the structure of sign group is given by $\imath(-1)=\omega$, $\delta(\omega)=0$, and
$\delta(t_i)=(i\; i+1)$, the permutation exchanging
$i$ and $i+1$ in $\set{1,\dots,n}$.

Below we use the homomorphisms
\begin{eqnarray*}
S^n\wedge-\colon\symt{m}\To\symt{n+m},\\
-\wedge S^m\colon\symt{n}\To\symt{n+m},
\end{eqnarray*}
defined on generators by
\begin{eqnarray*}
t_i\wedge S^m&=&t_i,\;\;1\leq i\leq n-1,\\
\omega\wedge S^m&=&\omega,\\
S^n\wedge t_i&=&t_{n+i},\;\;1\leq i\leq m-1,\\
S^n\wedge \omega&=&\omega.
\end{eqnarray*}

These homomorphisms are related by the following formula.

\begin{lem}\label{cambio}
Let $\tau_{n,m}\in\sym{n+m}$ be the permutation exchanging the first block of $n$ elements with the last
block of $m$ elements and let $\hat{\tau}_{n,m}\in\symt{n+m}$ be an element with
$\delta(\hat{\tau}_{n,m})=\tau_{n,m}$. Then for any $t\in \symt{n}$ we have 
\begin{eqnarray*}
(S^m\wedge t)\hat{\tau}_{n,m}&=&\hat{\tau}_{n,m}(t\wedge S^m)\omega^{nm\binom{\varepsilon\delta(t)}{2}}.
\end{eqnarray*}
\end{lem}

Notice that Lemma \ref{cambio} does not depend on the choice of $\hat{\tau}_{n,m}$ since the two possible
choices differ in $\omega$, which is central. For the proof of Lemma \ref{cambio} we choose
\begin{eqnarray}\label{tauhat}
\hat{\tau}_{n,m}&=&\overbrace{\!\!\!\underbrace{t_m\cdots t_1}_{m\text{  generators}}\!\!\!\cdots \cdots
\underbrace{t_{n+m-1}\cdots t_n}_{m\text{  generators}}}^{n\text{ groups of }m\text{  generators}}.
\end{eqnarray}

\begin{proof}[Proof of Lemma \ref{cambio}]
The equation holds for $t=\omega$, which is central, therefore we can restrict to the case $t=t_i$, $1\leq i\leq n-1$.
We check by induction in $j$ that
\begin{equation*}\tag{a}
t_{i+j}\cdots t_{i+1}t_it_{i+1}\cdots t_{i+j}\;\;=\;\;t_i\cdots t_{i+j-1}t_{i+j}t_{i+j-1}\cdots t_i.
\end{equation*}
For $j=1$ this is follows from (\ref{elados}), and if it is true for $j-1$
then
\begin{eqnarray*}
t_{i+j}t_{i+j-1}\cdots t_{i+1}t_it_{i+1}\cdots t_{i+j-1}t_{i+j}&=&t_{i+j}t_i\cdots t_{i+j-2}t_{i+j-1}t_{i+j-2}\cdots t_it_{i+j}\\
&=&t_i\cdots t_{i+j-2}t_{i+j}t_{i+j-1}t_{i+j}t_{i+j-2}\cdots t_i\\
\mbox{\scriptsize (\ref{elados})}&=&t_i\cdots t_{i+j-2}t_{i+j-1}t_{i+j}t_{i+j-1}t_{i+j-2}\cdots t_i.
\end{eqnarray*}
Equation (a) is equivalent to
\begin{equation*}\tag{b}
t_{i+j-1}\cdots t_it_{i+j}\cdots t_{i+1}t_i\;\;=\;\;t_{i+j}t_{i+j-1}\cdots t_{i}t_{i+j}\cdots t_{i+1}.
\end{equation*}
One can now easily check by using the other relations of the symmetric track groups that (b) for $j=m$ implies
\begin{eqnarray*}
\hat{\tau}_{n,m}t_i&=&t_{i+m}\hat{\tau}_{n,m}\omega^{m(n-2)},
\end{eqnarray*}
hence the lemma follows.
\end{proof}

\begin{rem}
The symmetric track groups were defined in \cite[5]{2hg2} in a geometric way in terms of tracks. In \cite[6]{2hg2} we relate them to the
positive pin group, obtaining in this way the presentation above, see \cite[Theorem 6.11]{2hg2}. 
The homomorphisms $S^n\wedge-$ and $-\wedge
S^m$ were geometrically defined in \cite[8]{2hg3}. We also give formulas for $S^n\wedge-$ 
in terms of the positive pin group in \cite[17]{2hg3}, from wich we derive the formulas for $S^n\wedge-$ in terms of the
presentation. The formulas for
$-\wedge S^m$ in terms of the presentation 
follow then from the definition of $-\wedge S^m$ in \cite[8]{2hg3} and from Lemma \ref{cambio}.
\end{rem}

%In the statement of this lemma $S^n\wedge-$ and $-\wedge
%S^m$ stand for the morphisms defined here in terms of the presentation.
%This lemma is satisfied by the morphisms $S^n\wedge-$ and $-\wedge
%S^m$ in \cite[8]{2hg3} from their very definition. As we remarked above the definition of $S^n\wedge -$ given
%here and in \cite[8]{2hg3} coincide, therefore the same follows for $-\wedge S^m$.

The next lemma encodes some relevant properties of the choices in (\ref{tauhat}).

\begin{lem}\label{recho}
The following equations hold for the elements in (\ref{tauhat}).
%\begin{eqnarray*}
%\hat{\tau}_{p,q}\hat{\tau}_{q,p}&=&\omega^{\binom{p}{2}\binom{q}{2}}.
%\end{eqnarray*}
\begin{enumerate}
\item $\hat{\tau}_{p,q}\hat{\tau}_{q,p}=\omega^{\binom{p}{2}\binom{q}{2}}$,
%\item $((S^r\wedge\hat{\tau}_{p,s})\wedge S^q)(S^{r+p}\wedge\hat{\tau}_{q,s})(\hat{\tau}_{p,r}\wedge
%S^{q+s})((S^p\wedge\hat{\tau}_{q,r})\wedge
%S^s)=$ 
%\newline $\hat{\tau}_{p+q,r+s}\omega^{rs\left(\binom{p}{2}+\binom{q}{2}+pq\right)}$.
\item $(S^r\wedge\hat{\tau}_{p,s}\wedge S^q)(S^{r+p}\wedge\hat{\tau}_{q,s})(\hat{\tau}_{p,r}\wedge
S^{q+s})(S^p\wedge\hat{\tau}_{q,r}\wedge
S^s)=$ 
\newline $\hat{\tau}_{p+q,r+s}\omega^{rs\left(\binom{p}{2}+\binom{q}{2}+pq\right)}$.

\end{enumerate}
\end{lem}

The proof only uses the presentation of the symmetric track groups as in Lemma \ref{cambio}. We leave it to the reader.

%\begin{lem}
%\begin{eqnarray*}
%(S^r\wedge\hat{\tau}_{p,s}\wedge S^q)(S^{r+p}\wedge\hat{\tau}_{q,s})(\hat{\tau}_{p,r}\wedge
%S^{q+s})(S^p\wedge\hat{\tau}_{q,r}\wedge S^s)&=&\hat{\tau}_{p+q,r+s}\omega^{rs\left(\binom{p}{2}+\binom{q}{2}+pq\right)}.
%\end{eqnarray*}
%\end{lem}

Relation (\ref{s5}) and Lemmas \ref{ader} and \ref{cambio} yield the following result.

\begin{lem}\label{tauel}
With the notation of Lemma \ref{cambio} the equation
\begin{eqnarray*}
[S^m\wedge t]\cdot[\tau_{n,m}]&
=
&
[\tau_{n,m}]\cdot[t\wedge S^m] %     +\binom{\varepsilon\delta(t)}{2}nmP(1|1)_H
\end{eqnarray*}
holds.
\end{lem}

We also use below the well-known cross product homomorphisms
\begin{equation*}
\sym{n}{\times}\sym{m}\To\sym{n+m}\colon (\sigma,\tau)\mapsto \sigma\times\tau.
\end{equation*}
Here $\sigma\times\tau$ permutes the first $n$ elements $\set{1,\dots,n}$ of 
$\set{1,\dots,n,n+1,\dots,n+m}$ according to $\sigma$ and the last $m$ elements $\set{n+1,\dots,n+m}$ according to
$\tau$. These homomorphisms satisfy $\tau_{n,m}(\sigma\times\tau)=(\tau\times\sigma)\tau_{n,m}$. Moreover, 
if $1_m\in\sym{m}$ denotes the unit of the symmetric group then $\delta(S^m\wedge t)=1_m\times\delta(t)$
and $\delta(t\wedge S^m)=\delta(t)\times 1_m$.

\section{$E_\infty$-quadratic pair algebras}
\renewcommand{\theequation}{E\arabic{equation}}\setcounter{equation}{0}

An \emph{$E_\infty$-quadratic pair algebra} is a 
quadratic pair algebra $B$ 
together with a cup-one product operation
$$\smile_1\colon B_{n,0}\times B_{m,0}\To B_{n+m,1},\;\;n,m\geq 0,$$
such that the quadratic pair module $B_{n,*}$ is a right $A(\symt{n})$-module
and the following compatibility conditions hold.
Let $x_i\in B_{n_i,0}$, $s_i\in B_{n_i,1}$, $a_i\in B_{n_i,ee}$, 
$g_i,g'_i\in\sym{n_i}$, and $r_i\in\symt{n_i}$. The product in the quadratic pair algebra $B$ is equivariant with
respect to the right $A(\symt{n})$-module structures in the following way
\begin{eqnarray}
\label{equi} (x_1\cdot[g_1])\cdot(x_2\cdot[g_2])&=&(x_1\cdot x_2)\cdot [g_1\times g_2],\\
\nonumber (s_1\cdot[g_1])\cdot(x_2\cdot[g_2])&=&(s_1\cdot x_2)\cdot [g_1\times g_2],\\
\nonumber (a_1\cdot([g_1]|[g_1'])_H)\cdot(a_2\cdot([g_2]|[g_2'])_H)&=&
\nonumber (a_1\cdot a_2)\cdot ([g_1\times g_2]| [g_1'\times g_2'])_H,\\
\nonumber x_1\cdot (x_2\cdot[r_2])&=&(x_1\cdot x_2)\cdot[S^{n_1}\wedge r_2].
\end{eqnarray}
The cup-one product measures the lack of commutativity, i.e. if $\tau_{p,q}\in\sym{p+q}$ denotes the permutation
exchanging the blocks $\set{1,\dots,p}$ and $\set{p+1,\dots,p+q}$, $p,q\geq 0$, then
\begin{eqnarray}
\label{lc}&&\\
\nonumber  (x_2\cdot x_1)\cdot[\tau_{n_1,n_2}]+\partial(x_1\smile_1x_2)&=&x_1\cdot x_2+\partial P(H(x_2)\cdot
TH(x_1))\cdot[\tau_{n_1,n_2}],\\
\nonumber (x_2\cdot s_1)\cdot[\tau_{n_1,n_2}]+\partial(s_1)\smile_1x_2&=&
s_1\cdot x_2+P(H(x_2)\cdot TH\partial(s_1))\cdot[\tau_{n_1,n_2}].
\end{eqnarray}
The cup-one product is itself commutative in the following sense
\begin{eqnarray}
\label{c1c} (x_2\smile_1x_1)\cdot[\tau_{n_1,n_2}]+x_1\smile_1x_2&=&-P(TH(x_1)\cdot H(x_2))\\
\nonumber &&+P(H(x_2)\cdot
TH(x_1))\cdot[\tau_{n_1,n_2}].
\end{eqnarray}
Let $1_n\in\sym{n}$ be the unit element.
The cup-one product also satisfies the following rules with respect to addition
\begin{eqnarray}
\label{c1+} x_1\smile_1(x_2+x_3)&=&x_1\smile_1x_2+x_1\smile_1x_3\\
\nonumber &&+P(\partial(x_1\smile_1 x_2)|(x_3\cdot x_1)\cdot[\tau_{n_1,n_3}])_H,%\\
\end{eqnarray}
multiplication
\begin{eqnarray}
\label{c1m} &&\\
\nonumber (x_1\cdot x_2)\smile_1x_3&=&((x_1\smile_1 x_3)\cdot x_2)\cdot[1_{n_1}\times\tau_{n_2,n_3}]+x_1\cdot (x_2\smile_1x_3)\\
\nonumber &&+P((\partial(x_1\smile_1 x_3)|(x_3\cdot x_1)\cdot[\tau_{n_1,n_3}])_H\cdot
H(x_2))\cdot[1_{n_1}\times\tau_{n_2,n_3}]\\
\nonumber &&+P(H(x_3)\cdot(x_1|x_1)_H\cdot TH(x_2))\cdot[\tau_{n_1+n_2,n_3}]\\
\nonumber &&-P((x_1|x_1)_H\cdot H(x_3)\cdot TH(x_2))\cdot [1_{n_1}\times\tau_{n_2,n_3}], 
\end{eqnarray}
and symmetric group action
\begin{eqnarray}
\label{c1e} (x_1\cdot[g_1])\smile_1(x_2\cdot[g_2])&=&(x_1\smile_1x_2)\cdot[g_1\times g_2].
\end{eqnarray}

If $B$ is an $E_\infty$-quadratic pair algebra then $h_0B$ is a commutative ring and $h_1B$ is an
$h_0B$-module, see \cite[Lemma 9.9]{2hg4}.

\numberwithin{equation}{section}
\renewcommand{\theequation}{\thesection.\arabic{equation}}
\setcounter{equation}{0}

Appart from Massey products the homology of an $E_\infty$-quadratic pair algebra is endowed with the following secondary
operation.

\begin{defn}\label{sq1}
Let $B$ be an $E_\infty$-quadratic pair algebra. Given an element $a\in h_{0}B_{2n,*}$ we define the \emph{cup-one square} of $a$
$$\sq_1(a)\in h_1B_{4n,*}$$
in the following way. 
Choose a representative $\bar{a}\in B_{2n,0}$ of $a$ and an element in the symmetric track group $\hat{\tau}\in \symt{4n}$ whose
boundary is the shuffle permutation
$\delta(\hat{\tau})=\tau_{2n,2n}$. 
Then $$\sq_1(a)=-\bar{a}^2\cdot[\hat{\tau}]+\bar{a}\smile_1\bar{a}-P(H(\bar{a})\cdot
TH(\bar{a}))\cdot[\tau_{2n,2n}]\in h_1B_{4n,*}.$$
\end{defn}

We leave it to the reader to check that the cup-one square does not depend on the choice of $\bar{a}$. However it does depend
on the choice of $\hat{\tau}$. There are two possible choices, namely $\hat{\tau}$ and
$\omega\hat{\tau}$. The difference between the two possible cup-one squares is computed in the following lemma.

\begin{lem}[{\cite[Lemma 9.11]{2hg4}}]\label{lor2}
Let $\sq_1$ be cup-one square in an $E_\infty$-quadratic pair algebra $B$ associated to the lift $\hat{\tau}$ of the shuffle
permutation and let $\sq_1^\omega$ be the cup-one square associated to $\omega\hat{\tau}$. Then given $a\in h_0B_{2n,*}$
\begin{eqnarray*}
\sq^\omega_1(a)&=&\sq_1(a)+a^2\cdot\eta.
\end{eqnarray*}
\end{lem}

The main result in \cite{2hg4} concerning $E_\infty$-quadratic pair algebras and ring spectra is the following.

\begin{thm}[{\cite[Theorem 9.12]{2hg4}}]\label{ca2}
There is a commutative diagram of functors
$$\xymatrix{{\left(\begin{array}{c}
\text{\emph{connective commutative}}\\
\text{\emph{ring spectra}}
\end{array}\right)}
\ar[d]_{\text{inclusion}}\ar[r]^-{\pi_{*,*}}&
{\left(E_\infty\text{\emph{-quadratic pair algebras}}\right)}
\ar[d]^{\text{forget}}\\
{\left(\text{\emph{connective ring spectra}}\right)}
\ar[r]^-{\pi_{*,*}}&
{\left(\text{\emph{quadratic pair algebras}}\right)}}$$
Here the lower arrow is the functor in Theorem \ref{a2}. Moreover, for a commutative ring spectrum $R$
the algebraic cup-one squares in $\pi_{*,*} R$ correspond to the topologically-defined cup-one
squares in $\pi_*R$.
\end{thm}

Theorem \ref{main3} follows from Theorem \ref{ca2} and from the following result.

\begin{thm}\label{moin3}
If $B$ is an $E_\infty$-quadratic pair algebra then the $k$-invariant (\ref{qpak}),
the Massey products in Definition \ref{qpamas}, and the cup-one squares in Definition \ref{sq1} associated to the
choices of $\hat{\tau}_{2n,2n}$ in (\ref{tauhat}) endow $h_0B$ with the structure of a commutative
ring with commutative secondary operations with coefficients in $h_1B$ in the sense of Definition \ref{crmp}.
\end{thm}

In order to prove this theorem we need a technical lemma.
%Notice that many of the equations above greatly simplify in case $H(x_i)=0$. This will be the case
%in many applications. 
$E_\infty$-quadratic pair algebras are defined above by using a minimal set of equations. Some other useful equations are listed
in the following lemma.

\begin{lem}\label{otros}
With the notation above the following equations are also satisfied in the $E_\infty$-quadratic pair algebra $B$
for elements with $H(x_i)=0$.
\begin{enumerate}
\item $(x_1\cdot[g_1])\cdot(s_2\cdot[g_2])=(x_1\cdot s_2)\cdot[g_1\times g_2]$,

\item $(x_1\cdot[r_1])\cdot x_2=(x_1\cdot x_2)\cdot[r_1\wedge
S^{n_2}]$,

\item $(s_2\cdot x_1)\cdot[\tau_{n_1,n_2}]+x_1\smile_1\partial(s_2)=x_1\cdot s_2$,

\item $(x_1+x_2)\smile_1x_3=x_1\smile_1x_3+x_2\smile_1x_3
+P(\partial(x_1\smile_1 x_3)|(x_3\cdot x_2)\cdot[\tau_{n_2,n_3}])_H$,

\item $x_1\smile_1(x_2\cdot x_3)= (x_2\cdot(x_1\smile_1x_3))\cdot[\tau_{n_1,n_2}\times 1_{n_3}]+(x_1\smile_1x_2)\cdot
x_3$,

\item $x_1\smile_1(x_2\cdot x_3)=((x_3\cdot x_1)\smile_1x_2)\cdot[\tau_{n_1+n_2,n_3}]+(x_1\cdot x_2)\smile_1x_3$.
\end{enumerate}
\end{lem}

\begin{proof}
The following equations hold.
$$\begin{array}{l}
((x_1\cdot [g_1])\cdot(s_2\cdot[g_2]))\cdot[\tau_{n_2,n_1}]\\
\mbox{\scriptsize (\ref{lc})}\qquad\qquad =
(s_2\cdot[g_2])\cdot(x_1\cdot [g_1])-\partial(s_2\cdot[g_2])\smile_1(x_1\cdot [g_1])\\
\qquad\qquad\qquad\quad +P(\underbrace{H(x_1\cdot [g_1])}_{\mbox{\scriptsize (\ref{hm},\ref{s1})}\;\;=0}\cdot TH\partial(s_2\cdot[g_2]))\cdot[\tau_{n_2,n_1}]\\
\mbox{\scriptsize (\ref{equi},\ref{dm},\ref{c1e})}\hspace{16pt} =(s_2\cdot x_1)\cdot[g_2\times g_1]-(\partial(s_2)\smile_1 x_1)\cdot[g_2\times g_1]\\
\mbox{\scriptsize (\ref{ldl},\ref{s1})}\hspace{29pt} =(s_2\cdot x_1-\partial(s_2)\smile_1x_1)\cdot[g_2\times g_1]\\
\mbox{\scriptsize (\ref{lc})}\hspace{40pt} =(x_1\cdot s_2)\cdot[\tau_{n_2,n_1}]\cdot[g_2\times g_1]\\
\mbox{\scriptsize (\ref{s3})}\hspace{41pt} =(x_1\cdot s_2)\cdot[g_1\times g_2]\cdot[\tau_{n_2,n_1}].
\end{array}$$
Now we obtain (1) multiplying by $[\tau_{n_1,n_2}]$ on the right and using (\ref{s3}).

Equation (2) follows from
\begin{eqnarray*}
(x_1\cdot x_2)\cdot[r_1\wedge S^{n_2}]&\st{\mbox{\scriptsize (\ref{lc})}}=&((x_2\cdot
x_1)\cdot[\tau_{n_1,n_2}]+\partial(x_1\smile_1x_2))\cdot[r_1\wedge S^{n_2}]\\
\mbox{\scriptsize (\ref{ldl},\ref{c1c})}&=& 
(x_2\cdot x_1)\cdot[\tau_{n_1,n_2}]\cdot[r_1\wedge S^{n_2}]+\partial(x_1\smile_1x_2)\cdot[r_1\wedge S^{n_2}]\\
&&+\overbrace{P((\partial(x_1\smile_1x_2)|(x_2\cdot x_1)\cdot[\tau_{n_1,n_2}])_H\cdot
H\partial[r_1\wedge S^{n_2}])}^{=\mbox{\scriptsize
(a)}}\\
\mbox{\scriptsize (\ref{tauel},\ref{equi},\ref{dm},\ref{s3},\ref{rl})}&=&
(x_2\cdot
(x_1\cdot[r_1]))\cdot[\tau_{n_1,n_2}]\hspace{-20pt}\underbrace{\pm}_{\mbox{\scriptsize according to
}\varepsilon\delta(r_1)}\hspace{-20pt}(x_1\smile_1x_2)\\
&&-(x_1\cdot[\delta(r_1)])\smile_1x_2+\mbox{(a)}\\
\mbox{\scriptsize (\ref{cem},\ref{otros}.4,\ref{dm},\ref{s3},\ref{pml},\ref{cem},\ref{hdt})}&=&(x_2\cdot
(x_1\cdot[r_1]))\cdot[\tau_{n_1,n_2}]+(x_1\cdot\partial[r_1])\smile_1x_2\\
%   &&+\binom{\varepsilon\delta(r_1)}{2}n_1n_2P(x_2\cdot x_1|x_2\cdot x_1)_H\\
\mbox{\scriptsize (\ref{lc},\ref{dm})}&=&(x_1\cdot[r_1])\cdot x_2.   %   +\binom{\varepsilon\delta(r_1)}{2}n_1n_2(x_1\cdot x_2\cdot\eta).
\end{eqnarray*}

Equation (3) follows from 
\begin{eqnarray*}
x_1\smile_1\partial(s_2)&\st{\mbox{\scriptsize (\ref{c1c})}}=&-(\partial(s_2)\smile_1 x_1)\cdot[\tau_{n_1,n_2}]\\
\mbox{\scriptsize (\ref{lc})}\quad&=&-(-(x_1\cdot s_2)\cdot[\tau_{n_2,n_1}]+s_2\cdot x_1)\cdot[\tau_{n_1,n_2}]\\
\mbox{\scriptsize (\ref{ldl}, \ref{s1}, \ref{s3})}\quad&=&
-(s_2\cdot x_1)\cdot[\tau_{n_1,n_2}]+x_1\cdot s_2.
\end{eqnarray*}

Equation (4) follows from
\begin{eqnarray*}
(x_1+x_2)\smile_1x_3&\st{\mbox{\scriptsize (\ref{lc})}}=&-(x_3\smile_1(x_1+x_2))\cdot[\tau_{n_1,n_3}]\\
\mbox{\scriptsize (\ref{c1+},\ref{ldl},\ref{s1})} &=&-(x_3\smile_1x_2)\cdot[\tau_{n_1,n_3}]
-(x_3\smile_1x_1)\cdot[\tau_{n_1,n_3}]\\
&&-P(\partial(x_3\smile_1x_1)|(x_2\cdot x_3)\cdot[\tau_{n_3,n_1}])_H\cdot[\tau_{n_1,n_3}]\\
\mbox{\scriptsize (\ref{c1c},\ref{pmr},\ref{cem},\ref{s1},\ref{s3})}&=&x_2\smile_1x_3+x_1\smile_1x_3\\
&&-P(\partial(x_3\smile_1x_1)\cdot[\tau_{n_1,n_3}]|x_2\cdot x_3)_H\\
\mbox{\scriptsize (\ref{qpm3})}&=&x_1\smile_1x_3+x_2\smile_1x_3\\
&&-P(\partial(x_1\smile_1x_3)|\!\!\!\!\!\!\!\!\!\!\!\!\!\!\!
\underbrace{\partial(x_2\smile_1x_3)}_{\mbox{\scriptsize (\ref{lc})}\;\;=-(x_3\cdot
x_2)\cdot[\tau_{n_1,n_3}]+x_2\cdot x_3}\!\!\!\!\!\!\!\!\!\!\!\!\!\!\!)_H\\
&&-P(\partial((x_3\smile_1x_1)\cdot[\tau_{n_1,n_3}])|x_2\cdot x_3)_H\\
\mbox{\scriptsize (\ref{c1c})}&=&x_1\smile_1x_3+x_2\smile_1x_3\\
&&+P(\partial(x_1\smile_1x_3)|x_3\cdot x_2\cdot[\tau_{n_2,n_3}])_H.
\end{eqnarray*}

Equation (5) follows from
\begin{eqnarray*}
x_1\smile_1(x_2\cdot x_3)&\st{\mbox{\scriptsize (\ref{c1c})}}=&-((x_2\cdot
x_3)\smile_1x_1)\cdot[\tau_{n_1,n_2+n_3}]\\
\mbox{\scriptsize (\ref{c1m},\ref{ldl},\ref{s1})} &=&-(x_2\cdot(x_3\smile_1 x_1))\cdot[\tau_{n_1,n_2+n_3}]\\
&&-((x_2\smile_1x_1)\cdot x_3)\cdot[1_{n_2}\times\tau_{n_3,n_1}]\cdot[\tau_{n_1,n_2+n_3}]\\
\mbox{\scriptsize (\ref{c1c},\ref{equi},\ref{rl},\ref{ldl},\ref{s1})} &=&
(x_2\cdot(x_1\smile_1 x_3))\cdot[1_{n_2}\times\tau_{n_3,n_1}]\cdot[\tau_{n_1,n_2+n_3}]\\
&&+((x_1\smile_1x_2)\cdot x_3)\cdot[\tau_{n_2,n_1}\times1_{n_3}]\cdot[1_{n_2}\times\tau_{n_3,n_1}]\cdot[\tau_{n_1,n_2+n_3}]\\
\mbox{\scriptsize (\ref{s3})}  &=&(x_2\cdot(x_1\smile_1x_3))\cdot[\tau_{n_1,n_2}\times 1_{n_3}]+(x_1\smile_1x_2)\cdot
x_3.
\end{eqnarray*}

Equation (6) follows from
$$\begin{array}{l}
((x_3\cdot x_1)\smile_1x_2)\cdot[\tau_{n_1+n_2,n_3}]+(x_1\cdot x_2)\smile_1x_3\\
\mbox{\scriptsize (\ref{c1m})}\quad =((x_3\smile_1x_2)\cdot
x_1)\cdot[1_{n_3}\times\tau_{n_1,n_2}]\cdot[\tau_{n_1+n_2,n_3}]+(x_3\cdot(x_1\smile_1x_2))\cdot[\tau_{n_1+n_2,n_3}]\\
\qquad\qquad +((x_1\smile _1x_3)\cdot x_2)\cdot[1_{n_1}\times\tau_{n_2,n_3}]+x_1\cdot(x_2\smile_1x_3)\\
\mbox{\scriptsize (\ref{otros}.5)}\; =-((x_2\smile_1x_3)\cdot x_1)\cdot[\tau_{n_1,n_2+n_3}]\\
\qquad\qquad +(x_1\smile_1(x_3\cdot x_2))\cdot[1_{n_1}\times\tau_{n_2,n_3}]+x_1\cdot(x_2\smile_1x_3)\\
\mbox{\scriptsize (\ref{qpm3})}\;\;\; =(x_1\smile_1(x_3\cdot x_2))\cdot[1_{n_1}\times\tau_{n_2,n_3}]+x_1\smile_1\partial(x_2\smile_1x_3)\\
\qquad\qquad +P(\partial(x_1\smile_1(x_3\cdot x_2))\cdot[1_{n_1}\times\tau_{n_2,n_3}]|
\partial((x_2\smile_1x_3)\cdot x_1)\cdot[\tau_{n_1,n_2+n_3}])_H\\
\mbox{\scriptsize (\ref{lc},\ref{c1+},\ref{c1e})} =x_1\smile_1(x_2\cdot x_3)
\end{array}$$
\end{proof}

We are now ready to prove Theorem \ref{moin3}.

\begin{proof}[Proof of Theorem \ref{moin3}]

We assume without loss of generality that all representatives chosen in $B_{*,0}$ are in $\ker H$.

(T7) By (\ref{lc}, \ref{dm}) we can take
\begin{eqnarray*}
\overline{ab}&=&\overline{ba}\cdot[\tau_{\abs{a},\abs{b}}]+\bar{a}\smile_1\bar{b},\\
\overline{bc}&=&\overline{cb}\cdot[\tau_{\abs{b},\abs{c}}]+\bar{b}\smile_1\bar{c},
\end{eqnarray*}
and so we do in this proof,
therefore
\begin{eqnarray*}
\mbox{(a)}\quad-\overline{ab}\cdot\bar{c}+\bar{a}\cdot\overline{bc}\!\!\!\!&\st{\mbox{\scriptsize
(\ref{lc},\ref{otros}.3)}}=&\!\!\!\!
-(\bar{a}\cdot\bar{b})\smile_1\bar{c}-(\bar{c}\cdot\overline{ab})\cdot[\tau_{\abs{a}+\abs{b},\abs{c}}]\\
&&\!\!\!\!+(\overline{bc}\cdot\bar{a})\cdot[\tau_{\abs{a},\abs{b}+\abs{c}}]+\bar{a}\smile_1(\bar{b}\cdot\bar{c})\\
&=&\!\!\!\!-(\bar{a}\cdot\bar{b})\smile_1\bar{c}-(\bar{c}\cdot(\overline{ba}\cdot[\tau_{\abs{a},\abs{b}}]+\bar{a}\smile_1\bar{b}))\cdot[\tau_{\abs{a}+\abs{b},\abs{c}}]\\
&&\!\!\!\!+((\overline{cb}\cdot[\tau_{\abs{b},\abs{c}}]+\bar{b}\smile_1\bar{c})\cdot\bar{a})\cdot[\tau_{\abs{a},\abs{b}+\abs{c}}]+\bar{a}\smile_1(\bar{b}\cdot\bar{c})\\
\mbox{(\scriptsize \ref{rl},\ref{ldl},\ref{equi},\ref{otros}.1)}\!\!\!\!\!\!\!&=&\!\!\!\!
-(\bar{a}\cdot\bar{b})\smile_1\bar{c}-(\bar{c}\cdot(\bar{a}\smile_1\bar{b}))\cdot[\tau_{\abs{a}+\abs{b},\abs{c}}]\\
&&\!\!\!\!+(-\bar{c}\cdot\overline{ba}+\overline{cb}\cdot\bar{a})\cdot[\tau_{\abs{b},\abs{c}}\times1_{\abs{a}}]\cdot[\tau_{\abs{a},\abs{b}+\abs{c}}]\\
&&\!\!\!\!+((\bar{b}\smile_1\bar{c})\cdot\bar{a})\cdot[\tau_{\abs{a},\abs{b}+\abs{c}}]+\bar{a}\smile_1(\bar{b}\cdot\bar{c})
\end{eqnarray*}
The element $-\bar{c}\cdot\overline{ba}+\overline{cb}\cdot\bar{a}$ represents $-\grupo{c,b,a}$, so it is in $\ker\partial$, in particular by Lemma \ref{mame} 
$$(-\bar{c}\cdot\overline{ba}+\overline{cb}\cdot\bar{a})\cdot[\tau_{\abs{b},\abs{c}}\times1_{\abs{a}}]\cdot[\tau_{\abs{a},\abs{b}+\abs{c}}]=(-1)^{\abs{a}\abs{b}+\abs{b}\abs{c}+\abs{c}\abs{a}+1}(-\overline{cb}\cdot\bar{a}+\bar{c}\cdot\overline{ba}).$$
Since $\ker\partial$ is central we only need to see that the rest of factors in the previous equation cancel, and
this follows from (\ref{otros}.6) since
$$\begin{array}{l}
-(\bar{c}\cdot(\bar{a}\smile_1\bar{b}))\cdot[\tau_{\abs{a}+\abs{b},\abs{c}}]
+((\bar{b}\smile_1\bar{c})\cdot\bar{a})\cdot[\tau_{\abs{a},\abs{b}+\abs{c}}]\\
\mbox{\scriptsize (\ref{c1c})}\qquad =-(\bar{c}\cdot(\bar{a}\smile_1\bar{b}))\cdot[\tau_{\abs{a}+\abs{b},\abs{c}}]
-((\bar{c}\smile_1\bar{b})\cdot\bar{a})\cdot[\tau_{\abs{b},\abs{c}}\times1_{\abs{a}}]\cdot[\tau_{\abs{a},\abs{b}+\abs{c}}]\\
\mbox{\scriptsize (\ref{c1m})}\qquad =-((\bar{c}\cdot\bar{a})\smile_1\bar{b})\cdot[\tau_{\abs{a}+\abs{b},\abs{c}}].
\end{array}$$

(T8) The following equation is obtained from the first equality in (a) above by inserting in the middle two elements which cancel
\begin{eqnarray*}
-\overline{ab}\cdot\bar{c}+\bar{a}\cdot\overline{bc}&=&
-(\bar{a}\cdot\bar{b})\smile_1\bar{c}-(\bar{c}\cdot\overline{ab})\cdot[\tau_{\abs{a}+\abs{b},\abs{c}}]+(\overline{ca}\cdot\bar{b})\cdot[\tau_{\abs{a}+\abs{b},\abs{c}}]\\
&&-\hspace{-46pt}\underbrace{(\overline{ca}\cdot\bar{b})}_{\mbox{\scriptsize
(\ref{lc})}\;=(\bar{b}\cdot\overline{ca})\cdot[\tau_{\abs{c}+\abs{a},\abs{b}}]+(\bar{c}\cdot\bar{a})\smile_1\bar{b}}\hspace{-44pt}\cdot[\tau_{\abs{a}+\abs{b},\abs{c}}]+(\overline{bc}\cdot\bar{a})\cdot[\tau_{\abs{a},\abs{b}+\abs{c}}]+\bar{a}\smile_1(\bar{b}\cdot\bar{c}).
\end{eqnarray*}
By Lemma \ref{mame}
\begin{eqnarray*}
(-\bar{c}\cdot\overline{ab}+\overline{ca}\cdot\bar{b})\cdot[\tau_{\abs{a}+\abs{b},\abs{c}}]&\in&
-(-1)^{\abs{a}\abs{c}+\abs{b}\abs{c}}\grupo{c,a,b},\\
(-\bar{b}\cdot\overline{ca}+\overline{bc}\cdot\bar{a})\cdot[\tau_{\abs{a},\abs{b}+\abs{c}}]&\in&
-(-1)^{\abs{a}\abs{b}+\abs{a}\abs{c}}\grupo{b,c,a},
\end{eqnarray*}
therefore, since $\ker\partial$ is central, (T8) follows from (\ref{otros}.6).

Let us now check simultaneously (T9) and (T10). By (\ref{lc}, \ref{dm}) we can take
\begin{eqnarray*}
%\overline{ab}&=&\overline{ba}\cdot[\tau_{\abs{a},\abs{b}}]+\bar{a}\smile_1\bar{b},\\
\overline{(2a)}&=&\bar{a}+\bar{a}.
\end{eqnarray*}
Moreover, by (\ref{s4},\ref{lc},\ref{dm},\ref{rl}) for $\abs{a}$ odd we can take
\begin{eqnarray*}
\overline{a(2a)}\;\;=\;\;\overline{(2a)a}&=&-\bar{a}^2\cdot[\hat{\tau}_{\abs{a},\abs{a}}]+\bar{a}\smile_1\bar{a}.
\end{eqnarray*}

\begin{eqnarray*}
-\overline{ab}\cdot \bar{a}+\bar{a}\cdot\overline{ba}\hspace{-15pt}&
=&\hspace{-15pt}-(\bar{a}\smile_1\bar{b})\cdot\bar{a}-(\overline{ba}\cdot\bar{a})\cdot[\tau_{\abs{a},\abs{b}}\times1_{\abs{a}}]\\
\mbox{\scriptsize (\ref{rl},\ref{equi},\ref{otros}.3)}\hspace{-15pt}&&\hspace{-15pt}+(\overline{ba}\cdot\bar{a})\cdot
\hspace{-35pt}\underbrace{[\tau_{\abs{a},\abs{b}+\abs{a}}]}_{\mbox{\scriptsize (\ref{s3})}\;=[1_{\abs{b}}\times{\tau}_{\abs{a},\abs{a}}]\cdot[\tau_{\abs{a},\abs{b}}\times1_{\abs{a}}]}
\hspace{-35pt}+\bar{a}\smile_1(\bar{b}\cdot\bar{a})\\
\mbox{\scriptsize
(\ref{rl},\ref{s4})}\hspace{-15pt}&=&\hspace{-15pt}
-(\bar{a}\smile_1\bar{b})\cdot\bar{a}-\mbox{(c)}\}=\left\{\begin{array}{l}
(\overline{ba}\cdot(\bar{a}+\bar{a}))
\cdot[\tau_{\abs{a},\abs{b}}\times1_{\abs{a}}],\,\mbox{$\abs{a}$ odd,}\\
0,\,\mbox{$\abs{a}$ even.}
\end{array}\right.\\
\hspace{-15pt}&&\hspace{-15pt}-(\overline{ba}\cdot\bar{a})\cdot\partial[S^{\abs{b}}\wedge\hat{\tau}_{\abs{a},\abs{a}}]\cdot[\tau_{\abs{a},\abs{b}}\times1_{\abs{a}}]
+\bar{a}\smile_1(\bar{b}\cdot\bar{a})\\
\mbox{\scriptsize
(\ref{dm},\ref{equi},\ref{c1m})}\hspace{-15pt}&=&\hspace{-15pt}-(\bar{a}\smile_1\bar{b})\cdot\bar{a}-\mbox{(c)}-(\bar{b}\cdot(\bar{a}^2\cdot[\hat{\tau}_{\abs{a},\abs{a}}]))\cdot[\tau_{\abs{a},\abs{b}}\times1_{\abs{a}}]\\
\hspace{-15pt}&&\hspace{-15pt}+(\bar{b}\cdot(\bar{a}\smile_1\bar{a}))\cdot[\tau_{\abs{a},\abs{b}}\times1_{\abs{a}}]+(\bar{a}\smile_1\bar{b})\cdot\bar{a}\\
\mbox{\scriptsize
(\ref{rl},\ref{ldl},\ref{s1})}\hspace{-15pt}&=&\hspace{-15pt}-(\bar{a}\smile_1\bar{b})\cdot\bar{a}-\mbox{(c)}\\
\hspace{-15pt}&&\hspace{-15pt}+(\bar{b}\cdot(-\bar{a}^2\cdot[\hat{\tau}_{\abs{a},\abs{a}}]+\bar{a}\smile_1\bar{a}))\cdot[\tau_{\abs{a},\abs{b}}\times1_{\abs{a}}]+(\bar{a}\smile_1\bar{b})\cdot\bar{a}
\\
\begin{array}{r}
\mbox{\scriptsize $\ker\partial$ central,}\\
\mbox{\scriptsize
(\ref{mame})}\end{array}\hspace{-5pt}&\left\{\begin{array}{l}
=\\
{}\\
\in
\end{array}\right.&\hspace{-15pt}
\begin{array}{l}
(-1)^{\abs{a}\abs{b}}b\cdot\sq_1(a), \quad \mbox{for $\abs{a}$ even,}\\
{}\\
(-1)^{\abs{a}\abs{b}}\grupo{b,a,2a},\quad \mbox{for $\abs{a}$ odd.}
\end{array}
\end{eqnarray*}

In order to check (T11) let $a,b\in h_0B_{2n,*}$, $\tau=\tau_{2n,2n}$, and $\hat{\tau}=\hat{\tau}_{2n,2n}$.
Using the ``bilinearity mod $P$'' of the product and the cup-one product we obtain
\begin{equation*}\tag{a}
\begin{array}{l}
-(\bar{a}+\bar{b})^2\cdot[\hat{\tau}]+(\bar{a}+\bar{b})\smile_1(\bar{a}+\bar{b})\\
\mbox{\scriptsize (\ref{rl},\ref{ldl},\ref{c1+},\ref{otros}.4)}\qquad
=-\bar{b}^2\cdot[\hat{\tau}]-(\bar{a}\cdot\bar{b})\cdot[\hat{\tau}]-(\bar{b}\cdot\bar{a})\cdot[\hat{\tau}]-\bar{a}^2\cdot[\hat{\tau}]-\mbox{(b)}\\
\;\;\quad\qquad\qquad\qquad\qquad
+\bar{a}\smile_1\bar{a}+\bar{b}\smile_1\bar{a}+\bar{a}\smile_1\bar{b}+\bar{b}\smile_1\bar{b}+\mbox{(c)}.
\end{array}
\end{equation*}
The central elements (b) and (c) are
\begin{eqnarray*}
\mbox{(b)} &=& P((\bar{b}\cdot\bar{a}+\bar{a}\cdot\bar{b}+\bar{b}^2|\bar{a}^2)_H\cdot H\partial[\hat{\tau}])
+P((\bar{a}\cdot\bar{b}+\bar{b}^2|\bar{b}\cdot\bar{a})_H\cdot H\partial[\hat{\tau}])\\
&&+P((\bar{a}\cdot\bar{b}+\bar{b}^2|\bar{a}\cdot\bar{b})_H\cdot H\partial[\hat{\tau}]),\\
\mbox{(c)} &=& P(\partial(\bar{a}\smile_1\bar{a})|(\bar{a}\cdot\bar{b})\cdot[\tau])_H 
 +P(\partial(\bar{a}\smile_1\bar{b})|\bar{b}^2\cdot[\tau])_H\\
&& +P(\partial((\bar{a}+\bar{b})\smile_1\bar{a})|(\bar{b}\cdot(\bar{a}+\bar{b}))\cdot[\tau])_H.
\end{eqnarray*}
In the middle of equation (a) we find the formula for $\sq_1(a)$ which is central in $B_{4n,1}$, so we can move
it to the end of the equation, as (b) and (c). Moreover, by \cite[Lemma 7.4]{2hg2} and (\ref{s5},\ref{s6}) we have
\begin{equation*}\tag{d}
n\,P(1|1)_H\;\;=\;\;[\tau]\cdot[\hat{\tau}]+[\hat{\tau}].
\end{equation*}
This formula is used in the following equation.
\begin{eqnarray*}
\bar{b}\smile_1\bar{a}+\bar{a}\smile_1\bar{b}&\st{\mbox{\scriptsize
(\ref{c1c})}}=&-(\bar{a}\smile_1\bar{b})\cdot[\tau]+\bar{a}\smile_1\bar{b}\\
\mbox{\scriptsize (\ref{rl})}&=&(\bar{a}\smile_1\bar{b})\cdot(-[\tau]+1)\\
\mbox{\scriptsize (\ref{s4})}&=&(\bar{a}\smile_1\bar{b})\cdot\partial[\hat{\tau}]\\
\mbox{\scriptsize
(\ref{dm},\ref{lc})}&=&(-(\bar{b}\cdot\bar{a})\cdot[\tau]+\bar{a}\cdot\bar{b})\cdot[\hat{\tau}]\\
\mbox{\scriptsize
(\ref{ldl})}&=&-(\bar{b}\cdot\bar{a})\cdot([\tau]\cdot[\hat{\tau}])+(\bar{a}\cdot\bar{b})\cdot[\hat{\tau}]\\
&&+\underbrace{P((-\bar{a}\cdot\bar{b}+(\bar{b}\cdot\bar{a})\cdot[\tau]|(\bar{b}\cdot\bar{a})\cdot[\tau])_H\cdot
H\partial[\hat{\tau}])}_{= \mbox{\scriptsize (e)}}\\
\mbox{\scriptsize
(d,\ref{pml},\ref{cem})}&=&
(\bar{b}\cdot\bar{a})\cdot[\hat{\tau}]+(\bar{a}\cdot\bar{b})\cdot[\hat{\tau}]+\mbox{(e)}+n
\underbrace{P(\bar{b}\cdot\bar{a}|\bar{b}\cdot\bar{a})_H}_{=a\cdot b\cdot\eta}.
\end{eqnarray*}
This shows that (a) simplifies to give the following equation
\begin{eqnarray*}
\sq_1(a+b)&=&\sq_1(a)+\sq_1(b)+n\cdot a\cdot b\cdot\eta+\mbox{(b)}+\mbox{(c)}+\mbox{(e)}.
\end{eqnarray*}
Now one uses
(\ref{cem},\ref{s3},\ref{hdt}) and the elementary properties of quadratic pair modules to check that
\begin{eqnarray*}
\mbox{(b)}+\mbox{(c)}+\mbox{(e)}&=&P(\bar{b}\cdot\bar{a}|\bar{b}\cdot\bar{a})_H\;\;=\;\;a\cdot b\cdot\eta,
\end{eqnarray*}
hence we are done.

Finally (T12) will follow from (\ref{recho},\ref{lor2}) once we check that for the cup-one
square $\overline{\sq}_1(a\cdot b)$ associated to the following lift of $\tau_{\abs{a}+\abs{b},\abs{a}+\abs{b}}$
\begin{eqnarray*}
\tilde{\tau}&=&(S^{\abs{a}}\wedge\hat{\tau}_{\abs{b},\abs{a}}^{-1}\wedge S^{\abs{b}})(S^{2\abs{a}}\wedge\hat{\tau}_{\abs{b},\abs{b}})
(\hat{\tau}_{\abs{a},\abs{a}}\wedge
S^{2\abs{b}})(S^{\abs{a}}\wedge\hat{\tau}_{\abs{b},\abs{a}}\wedge S^{\abs{b}})
\end{eqnarray*}
the following formula holds.
\begin{eqnarray*}
\overline{\sq}_1(a\cdot b)&=&a^2\cdot\sq_1(b)+\sq_1(a)\cdot b^2.
\end{eqnarray*}

%In order to check this equation and 
%for the sake of simplicity we simply denote $\tau=\tau_{2n,2n}$ and $\hat{\tau}=\hat{\tau}_{2n,2n}$.
The following equation holds.
\begin{eqnarray*}
[\tilde{\tau}]\hspace{-10pt}&\st{\mbox{\scriptsize (\ref{s5})}}=&\hspace{-10pt}
[1_{\abs{a}}\times\tau_{\abs{a},\abs{b}} \times1_{\abs{b}}]\cdot[(S^{2\abs{a}}\wedge\hat{\tau}_{\abs{b},\abs{b}} )
(\hat{\tau}_{\abs{a},\abs{a}} \wedge
S^{2\abs{b}})(S^{\abs{a}}\wedge\hat{\tau}_{\abs{b},\abs{a}} \wedge S^{\abs{b}})]\\
&&+[S^{\abs{a}}\wedge\hat{\tau} ^{-1}_{\abs{b},\abs{a}}\wedge S^{\abs{b}}]\\
\mbox{\scriptsize (\ref{ader},\ref{rl})}\hspace{-10pt}&=&
\hspace{-10pt}[1_{\abs{a}}\times\tau_{\abs{a},\abs{b}}
\times1_{\abs{b}}]\cdot[(S^{2\abs{b}}\wedge\hat{\tau}_{\abs{b},\abs{b}} )
(\hat{\tau}_{\abs{a},\abs{a}} \wedge S^{2\abs{b}})]\cdot[1_{\abs{a}}\times\tau_{\abs{b},\abs{a}} \times1_{\abs{b}}]\\
&&\hspace{-10pt}+\underbrace{[1_{\abs{a}}\times\tau_{\abs{a},\abs{b}}
\times1_{\abs{b}}]\cdot[S^{\abs{a}}\wedge\hat{\tau}_{\abs{b},\abs{a}} \wedge
S^{\abs{b}}]+[S^{\abs{a}}\wedge\hat{\tau}^{-1}_{\abs{b},\abs{a}}\wedge S^{\abs{b}}]}_{\mbox{\scriptsize (\ref{10},\ref{s5})}\quad=0}\\
\mbox{\scriptsize
(\ref{s5})}&=&\hspace{-10pt}[1_{\abs{a}}\times\tau_{\abs{a},\abs{b}}
\times1_{\abs{b}}]\cdot([1_{2\abs{a}}\times\tau_{\abs{b},\abs{b}} ]\cdot[\hat{\tau}_{\abs{a},\abs{a}} \wedge
S^{2\abs{b}}]\\
&&+[S^{2\abs{a}}\wedge\hat{\tau}_{\abs{b},\abs{b}} ])\cdot[1_{\abs{a}}\times\tau_{\abs{b},\abs{a}}
\times1_{\abs{b}}].
\end{eqnarray*}
This equation is used below. We will also use the notation
\begin{eqnarray*}
\mbox{(*)}&=&P((\bar{a}\cdot\partial(\bar{b}\smile_1\bar{a})\cdot\bar{b}|(\bar{a}^2\cdot\bar{b}^2)\cdot[1_{\abs{a}}\times\tau_{\abs{b},\abs{a}}\times1_{\abs{b}}])_H\cdot H\partial[\tilde{\tau}])\\
\mbox{\scriptsize
(\ref{qpm3},\ref{hdt},\ref{cem})}&=&[(\bar{a}\cdot(\bar{b}\smile_1\bar{a})\cdot\bar{b})\cdot[\tau_{\abs{a}+\abs{b},\abs{a}+\abs{b}}],-(\bar{a}^2\cdot\bar{b}^2)\cdot[1_{\abs{a}}\times\tau_{\abs{b},\abs{a}}\times1_{\abs{b}}]\cdot[\tilde{\tau}]].
\end{eqnarray*}
Now the formula follows from the following equations.
$$\begin{array}{l}
-(\bar{a}\cdot \bar{b})^2\cdot[\tilde{\tau}]+(\bar{a}\cdot\bar{b})\smile_1(\bar{a}\cdot\bar{b})\\
\left(\begin{array}{c}
\mbox{\scriptsize \ref{lc},\ref{rl},\ref{ldl},\ref{equi},}\\
\mbox{\scriptsize\ref{c1m},\ref{otros}.4,\ref{hm},\ref{c1e}}
\end{array}\right)\quad
=-((\bar{a}^2\cdot\bar{b}^2)\cdot[1_{\abs{a}}\times\tau_{\abs{b},\abs{a}}\times1_{\abs{b}}]+\bar{a}\cdot\partial(\bar{b}\smile_1\bar{a})\cdot\bar{b})\cdot[\tilde{\tau}]\\
\qquad\qquad\qquad\qquad\quad\;\;
\left.\begin{array}{l}
+\underbrace{(\bar{a}\cdot(\bar{a}\smile_1\bar{b})\cdot\bar{b})\cdot
\overbrace{[\tau_{\abs{a},\abs{a}}\times1_{2\abs{b}}]\cdot[1_{\abs{a}}\times\tau_{\abs{b},\abs{a}+\abs{b}}]}^{\mbox{\scriptsize
(\ref{s3})}\;\;=[1_{\abs{a}}\times\tau_{\abs{b},\abs{a}}\times1_{\abs{b}}]\cdot[\tau_{\abs{a}+\abs{b},\abs{a}+\abs{b}}]}}_{=\mbox{\scriptsize(a)}}\\
+((\bar{a}\smile_1\bar{a})\cdot\bar{b}^2)\cdot\hspace{-30pt}\underbrace{[1_{\abs{a}}\times\tau_{\abs{b},\abs{a}+\abs{b}}]}_{\mbox{\scriptsize
(\ref{s3})}\;\;=[1_{2\abs{a}}\times\tau_{\abs{b},\abs{b}}]\cdot[1_{\abs{a}}\times\tau_{\abs{b},\abs{a}}\times1_{\abs{b}}]}\\
+(\bar{a}^2\cdot(\bar{b}\smile_1\bar{b}))\cdot[1_{\abs{a}}\times\tau_{\abs{b},\abs{a}}\times1_{\abs{b}}]\\
+\underbrace{\bar{a}\cdot(\bar{b}\smile_1\bar{a})\cdot\bar{b}}_{=\mbox{\scriptsize (b)}}
\end{array}\right\}=\mbox{(c)}\\
\mbox{\scriptsize (\ref{ldl},\ref{hdt})}\qquad\qquad\qquad
=-\hspace{-40pt}\underbrace{(\bar{a}\cdot\partial(\bar{b}\smile_1\bar{a})\cdot\bar{b})\cdot[\tilde{\tau}]}_{\begin{array}{l}\scriptstyle\mbox{\tiny
(\ref{dm})}\quad\;\;=(\bar{a}\cdot(\bar{b}\smile_1\bar{a})\cdot\bar{b})\cdot\partial[\tilde{\tau}]\\
\scriptstyle\mbox{\tiny
(\ref{s4},\ref{rl})}\;\;=-\hspace{-20pt}\underbrace{\scriptstyle(\bar{a}\cdot(\bar{b}\smile_1\bar{a})\cdot\bar{b})\cdot[\tau_{\abs{a}+\abs{b},\abs{a}+\abs{b}}]}_{\mbox{\tiny
cancels with (a) and (*) by (\ref{equi},\ref{c1c},\ref{rl},\ref{ldl})}}\\
\scriptstyle\qquad\quad\;\;+\underbrace{\scriptstyle(\bar{a}\cdot(\bar{b}\smile_1\bar{a})\cdot\bar{b})}_{\tiny\mbox{\tiny
cancels with (b)}}
\end{array}}\hspace{-40pt}
-(\bar{a}^2\cdot\bar{b}^2)\cdot[1_{\abs{a}}\times\tau_{\abs{b},\abs{a}}\times1_{\abs{b}}]\cdot[\tilde{\tau}]\\
\hspace{108pt}-\mbox{(*)}+\mbox{(c)}\\
\mbox{\scriptsize (\ref{equi},\ref{otros}.2,\ref{rl},\ref{ldl})}\qquad\quad
=(-\bar{a}^2\cdot(\bar{b}^2\cdot[\hat{\tau}_{\abs{b},\abs{b}}])
-(\bar{a}^2\cdot[\hat{\tau}_{\abs{a},\abs{a}}])\cdot(\bar{b}^2\cdot[\tau_{\abs{b},\abs{b}}])\\
\hspace{100pt}+(\bar{a}\smile_1\bar{a})\cdot(\bar{b}^2\cdot[\tau_{\abs{b},\abs{b}}])+\bar{a}^2\cdot(\bar{b}\smile_1\bar{b}))\cdot[1_{\abs{a}}\times_{\abs{b},\abs{a}}\tau\times1_{\abs{b}}]\\
\mbox{\scriptsize (\ref{hm},\ref{rl},\ref{ldl})}\qquad\qquad\;\;\;
=(\bar{a}^2\cdot(-\bar{b}^2\cdot[\hat{\tau}_{\abs{b},\abs{b}}]+\bar{b}\smile_1\bar{b})\\
\hspace{100pt}+(-\bar{a}^2\cdot[\hat{\tau}_{\abs{a},\abs{a}}]+\bar{a}\smile_1\bar{a})\cdot(\bar{b}^2\cdot[\tau_{\abs{b},\abs{b}}]))\cdot[1_{\abs{a}}\times\tau_{\abs{b},\abs{a}}\times1_{\abs{b}}]\\
\mbox{\scriptsize (\ref{mame})}\hspace{72pt}=a^2\cdot\sq_1(b)+\sq_1(a)\cdot b^2.
\end{array}$$
\end{proof}

% ----------------------------------------------------------------
\bibliographystyle{amsalpha}
\bibliography{Fernando}
\end{document}